\documentclass[final]{siamart171218}

\usepackage{amsfonts}
\usepackage{graphicx}
\usepackage{epstopdf}
\usepackage{algorithmic}
\usepackage[utf8]{inputenc}
\usepackage{amsmath}
\usepackage{amsfonts}
\usepackage{amssymb}
\usepackage{paralist}
\usepackage{gnuplot-lua-tikz}
\usepackage{tikz} 
    \usetikzlibrary{backgrounds}
    \usetikzlibrary{external}
\usepackage{mathrsfs}
\usepackage{xfrac}
\usepackage{commath}
\usepackage{array}
\usepackage{booktabs}
\usepackage{dcolumn}
\usepackage{afterpage}
\usepackage{tcolorbox}
\usepackage[super]{nth}
\ifpdf
  \DeclareGraphicsExtensions{.eps,.pdf,.png,.jpg}
\else
  \DeclareGraphicsExtensions{.eps}
\fi

\newcommand{\TheTitle}{End-Point Corrections for the Midpoint Rule}
\newcommand{\TheAuthors}{Ricardo L. U. F. Pinto, and Bernardo B. Monteiro}

\headers{\TheTitle}{\TheAuthors}

\title{{\TheTitle}\thanks{Submitted to the editors \today.}}

\author{
  Ricardo L. U. F. Pinto\thanks{Department of Mechanical Engineering, Federal University of Minas Gerais, Belo Horizonte, MG, Brazil
      (\email{utsch@demec.ufmg.br}, \email{b.b.monteiro@gmail.com}).}
  \and
  Bernardo B. Monteiro\footnotemark[2]
}

\newcommand{\R}{\mathbb{R}}
\newcommand{\N}{\mathbb{N}}
\newtcbox{\todo}[1][yellow]{on line,arc=7pt,colback=#1!10!white,colframe=#1!50!black,before upper={\rule[-3pt]{0pt}{10pt}},boxrule=1pt, boxsep=0pt,left=6pt,right=6pt,top=2pt,bottom=2pt}

\ifpdf
\hypersetup{
  pdftitle={\TheTitle},
  pdfauthor={\TheAuthors}
}
\fi

\begin{document}

\maketitle

% REQUIRED
\begin{abstract}
  In this paper we present a new family of rules for numerical integration.
  This family has up to half the error of the widely used Newton-Cotes rules
   when a sufficient number of points is evaluated
   and also much better numerical stability for high orders. 
  These rules can be written as the midpoint rule with a correction term, 
   providing a straightforward and computationally cheap way to obtain error estimations.
  The rules are interpolatory and use evenly spaced points,
   which makes them well suited for many practical applications.
  Their major potential disadvantage is the use of points outside the integration interval.

%  Neste trabalho propõe-se uma nova família de fórmulas
%  interpolatórias de integração numérica que apresenta acurácia
%  superior às fórmulas de Newton-Cotes, muito usadas na prática. Em
%  particular serão apresentadas formulas compostas equivalentes às
%  regras de Simpson e de Boole (3ª e 5ª ordem). As fórmulas utilizam
%  pontos externos ao intervalo de integração mas não requerem
%  informações sobre derivadas. Elas permitem estimar o erro de
%  integração a partir de poucas avaliações do integrando na vizinhança
%  das extemidades do intervalo de integração.
\end{abstract}

% REQUIRED

% REQUIRED

\section{Introduction}
Numerical integration is one of the most basic procedures used when tackling
practical problems in science, technology, engineering and mathematics.
Although there are many well established techniques for numerical
integration\cite{DavisRabinowitz}, some improvements have been suggested
throughout the years\cite{favati1991interpolatory,favati1995new,hale2008new,milovanovic2016generalized}.

A straightforward way to improve the accuracy of existing quadrature rules is
the use of endpoint corrections\cite{kapur1997high,rokhlin1990end}. Much
research has been done on deriving end-point corrections for the
trapezoidal~rule, but we feel that the midpoint~rule has been somewhat
neglected. It has the same order of accuracy and is arguably slightly more
precise than the former for polynomial behaving functions. Indeed, the error of the
midpoint~rule for $\int_a^bf(x)dx$ is given by $-\frac{(b-a)^3}{24N^2}f''(\xi)$,
while the trapezoidal~rule has an error given by
$-\frac{(b-a)^3}{12N^2}f''(\xi), \quad \xi \in [a,b]$.

In this paper, we present a new family of integration rules, which can be
written in the form of end-point corrections of arbitrary order for the
midpoint rule. These rules were derived with three main principles in mind:
\begin{enumerate}
\item the use of interpolatory polynomials;
\item the use of points outside the integration interval;
\item the overlapping of these points, yielding end-point like correction terms.
\end{enumerate}

We think these rules compare favorably to Newton-Cotes' rules. They have
less error when a large number of function evaluations are used and better numerical stability.
Furthermore, as they can
be presented as correction formulae, the correction term yields a good
heuristic for \emph{a priori} error estimation, with only very few function
evaluations.

% Outline
This paper is organized as follows. 
In \cref{sec:new_quadrature_rules}, 
 we derive the new family of integration rules, 
 starting with a didactic and intuitive derivation of the three point rule 
 in order to convey the spirit of the new rules to the reader. 
The derivation is then generalized to a rule with an arbitrary number of points, which is then particularized to yield a five point rule, 
 as we feel this result could be most useful for practical applications. 
In \cref{sec:integration_error}, we derive an error formula based on 
 Peano's theory\footnote{
     the reader should refer to 
     \cite{DavisRabinowitz} for an overview of the method and 
     \cite{favati1995peano} for some interesting results concerning 
      symmetric rules
 } for the general rule. 
The formula thus derived proves that the $n$ point rule is also of $n$th order,
 hence justifying the use of the expression ``$n$th order rule'' 
 instead of ``$n$ point rule'' in section titles.
In \cref{sec:modified_rules}, we tackle the rules' limitation of 
 using points outside the integration interval 
 by providing two alternatives for the third order case. 
Only one of them requires the use of derivatives, 
 and even then they need only be evaluated at the endpoints.
\Cref{sec:results} analyses the main results of this paper, 
 including a comparison with the widely used Simpson's rule and 
 heuristics for error estimation.
Finally, \cref{sec:conclusion} presents the conclusions of this paper.
We also included several numerical examples in \cref{sec:numerical_results}.

\section{New quadrature rules}
\label{sec:new_quadrature_rules}
\subsection{Derivation of a third order rule}
\label{sec:3rd_order}

\begin{figure}
	\centering
	\begin{tikzpicture}
    [scale = 3,
        axis/.style = {thin}
    ]
    \pgfmathsetmacro{\ya}{.5}  %y(-1)
    \pgfmathsetmacro{\yo}{1.2} %y(0)
    \pgfmathsetmacro{\yb}{1.3} %y(1)
    \pgfmathsetmacro{\yd}{1.28} %p_2(.5)
    \pgfmathsetmacro{\xmin}{-1.2}
    \pgfmathsetmacro{\xmax}{1.2}
    \pgfmathsetmacro{\ymin}{-.2}
    \pgfmathsetmacro{\ymax}{1.5}

    %f(x): 3rd degree polynomial
    \draw[domain=\xmin:\xmax][samples=100][very thick]
        plot (\x,{-(\x-1)*\x*(\x-.5)/3*\ya+ 2*(\x-1)*(\x+1)*(\x-.5)*\yo + \x*(\x+1)*(\x-.5)*\yb + -(\x-1)*\x*(\x+1)*4/1.5*\yd}) 
	[anchor=south west] node {$f(x)$};
    
    %Interpolating parabola
    \draw[domain=\xmin:\xmax][samples=100][semithick]
        plot (\x,{(\x-1)*\x/2*\ya+ -(\x-1)*(\x+1)*\yo + \x*(\x+1)/2*\yb}) 
	[anchor=north west] node {$p_2(x)$};

    %Interpolation points
    \foreach \p in {(-1,\ya), (0,\yo), (1,\yb)}
    \fill \p circle [radius=.025];
    
    \begin{scope}[on background layer]
        %Area corresponding to integration
        \fill[black!20!white][domain=-.5:.5]
             plot (\x,{(\x-1)*\x/2*\ya+ -(\x-1)*(\x+1)*\yo + \x*(\x+1)/2*\yb}) -- (.5,0) -- (-.5,0) --cycle;

	%Borders from the interval of integration
        \draw[help lines]
	    (-.5,0) -- (-.5,\ymax) (.5,0) -- (.5,\ymax);
    \end{scope}

    %Axis
	\draw[axis] (\xmin,0) -- (\xmax,0); %x axis
    %\foreach \x in {-1, -1/2, 1/2, 1}
    %    \draw[axis] (\x,.05) -- (\x,-.05) %x axis ticks
	%    [below] node {$\x$}; %x axis labels
	\draw[axis] (-1,.05) -- (-1,-.05)
		[below] node {$-h$}
		[axis] (-.5,.05) -- (-.5,-.05)
		[below] node{$-h/2$}
		[axis] (.5,.05) -- (.5,-.05)
		[below] node{$h/2$}
		[axis] (1,.05) -- (1,-.05)
		[below] node{$h$};
	\draw[axis] (0, \ymin) -- (0, \ymax); %y axis

\end{tikzpicture}
	
	\caption{Third order quadrature rule. $f(x)$ is approximated by $p_2(x)$,
    which is the interpolating polynomial constructed from points $-h, 0, h$.
    $p_2(x)$ is integrated from $-h/2$ to $h/2$ in order to approximate the
    integral of $f(x)$ over $[-h/2, h/2]$.} 
	\label{fig:3rd_order_simple}

\end{figure}
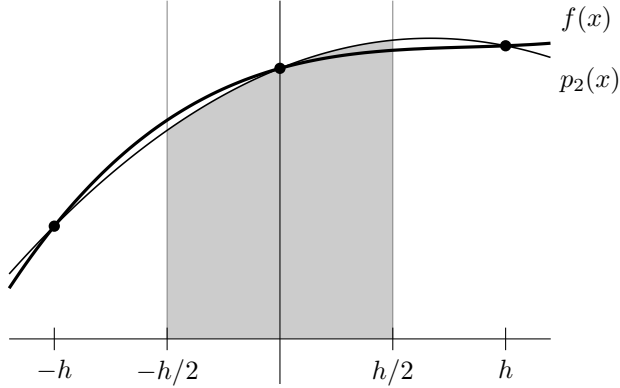

We can approximate the function $f$ by constructing a polynomial from a central
point $f(0)$ and two symmetric side points $f(-h)$ and $f(h)$, $h \in \R$, as show in
\cref{fig:3rd_order_simple}. Using Lagrange's interpolation formula \cite{Atkinson}, we get

\begin{equation}
p_2(x) = \frac{x(x-h)}{2h^2} f(-h) + \frac{(x-h)(x+h)}{-h^2} f(0)
+ \frac{x(x+h)}{2h^2} f(h)
\end{equation}

It is easy to check that $p_2(x)$ and $f(x)$ assume the same values in all three
interpolation points.
If we assume $p_2$ as an approximation of $f(x)$,
 the integral over the small interval $[-h/2,h/2]$ can then be written as

\begin{equation}
\int_{-h/2}^{h/2}f(x) \dif x \approx \int_{-h/2}^{h/2}p_2(x) \dif x
= h\frac{f(-h)-22f(0)+f(h)}{24}
\end{equation}

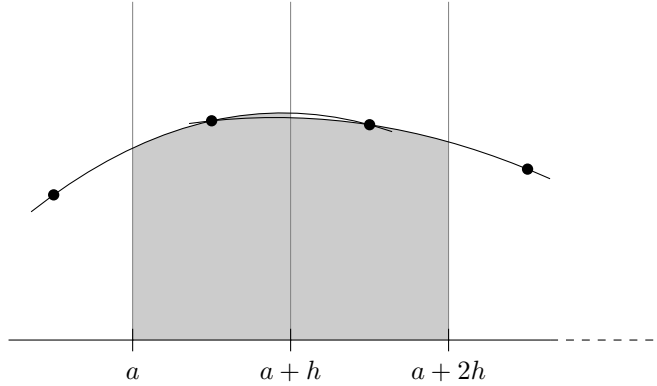
\begin{figure}
    \centering
    \begin{tikzpicture}[
        declare function = {f(\x) = .3*(sin((1.5*\x+.5) r)+1)+pow(\x+1,.3)-0.8;
                        p2(\x,\a,\ya,\b,\yb,\c,\yc) = %interpolating parabola
                        (\x-\b)*(\x-\c)/((\a-\b)*(\a-\c))*\ya +
                        (\x-\a)*(\x-\b)/((\c-\a)*(\c-\b))*\yc +
                        (\x-\a)*(\x-\c)/((\b-\a)*(\b-\c))*\yb ;
    },
    scale = 3,
    axis/.style = {thin}
    ]
    \pgfmathsetmacro{\xmin}{0}
    \pgfmathsetmacro{\xmax}{2.5}
    \pgfmathsetmacro{\epsilon}{.2} %extra extension for axis, etc
    \pgfmathsetmacro{\h}{.7} %integration step
    \pgfmathsetmacro{\ymax}{1.5}

    %f(x)
    %\draw[domain=\xmin-\epsilon:\xmin+3*\h+\epsilon/2][samples=40][very thick]
    %    plot (\x,{f(\x)});
    %\draw[domain=\xmin+3*\h+\epsilon/2:\xmax-2*\h-\epsilon/2][samples=20][very
    %thick, dotted]
    %    plot (\x,{f(\x)});
    %\draw[domain=\xmax-2*\h-\epsilon/2:\xmax+\epsilon][samples=40][very thick]
    %    plot (\x,{f(\x)})
	%[anchor=west] node {$f(x)$};

    %interpolationg parabolas
    \draw[domain=\xmin-\epsilon/2:\xmin+2*\h+\epsilon/2][samples=20][thin]
        plot
        (\x,{p2(\x,\xmin,f(\xmin),\xmin+\h,f(\xmin+\h),\xmin+2*\h,f(\xmin+2*\h)});
    \draw[domain=\xmin+\h-\epsilon/2:\xmin+3*\h+\epsilon/2][samples=20][thin]
        plot
        (\x,{p2(\x,\xmin+\h,f(\xmin+\h),\xmin+2*\h,f(\xmin+2*\h),\xmin+3*\h,f(\xmin+3*\h)});
    \draw[domain=\xmax-2*\h-\epsilon/2:\xmax+\epsilon/2][samples=20][thin];

    %Interpolation points
    \foreach \x in {\xmin,\xmin+\h,\xmin+2*\h,\xmin+3*\h}
    \fill (\x,{f(\x)}) circle [radius=.025];
    
    \begin{scope}[on background layer]
        %Area corresponding to integration
        \fill[black!20!white]
        [domain=\xmin+.5*\h:\xmin+1.5*\h] plot 
             (\x,{p2(\x,\xmin,f(\xmin),\xmin+\h,f(\xmin+\h),\xmin+2*\h,f(\xmin+2*\h))}) 
             -- (\xmin+1.5*\h,0) -- (\xmin+.5*\h,0) -- cycle
        [domain=\xmin+1.5*\h:\xmin+2.5*\h] plot 
             (\x,{p2(\x,\xmin+\h,f(\xmin+\h),\xmin+2*\h,f(\xmin+2*\h),\xmin+3*\h,f(\xmin+3*\h))}) 
             -- (\xmin+2.5*\h,0) -- (\xmin+1.5*\h,0) -- cycle;
        
        %Border from each step
        \foreach \x in {\xmin+.5*\h,\xmin+1.5*\h,\xmin+2.5*\h}
        \draw[help lines] (\x,0)--(\x,\ymax);
    \end{scope}

    %Axis
	\draw[axis] 
    (\xmin-\epsilon,0) -- (\xmin+3*\h+\epsilon/2,0) %x axis
	edge[dashed] (\xmax+\epsilon,0)
    %tick: a
    (\xmin+.5*\h,.05) -- (\xmin+.5*\h,-.05) 
            [below] node {$a\vphantom{b}$}
    %tick: a+h
    (\xmin+1.5*\h,.05) -- (\xmin+1.5*\h,-.05) 
            [below] node {$a+h$}
    %tick: a+2*h
    (\xmin+2.5*\h,.05) -- (\xmin+2.5*\h,-.05) 
            [below] node {$a+2h$};

\end{tikzpicture}
    \caption{Third order composite quadrature rule.
     The function is approximated by a polynomial 
     as shown in \cref{fig:3rd_order_simple} for each sub-interval.}
    \label{fig:3rd_order_composite}
\end{figure}

Using this formula in all sub-intervals of the form $[a+ih, a+(i+1)h]$,
$i=0,1,\dots,N-1$,  we get the composite rule 
\begin{multline}
\label{eqn:3rd_order_composite}
    \int_a^b f(x) \dif x \approx h\frac{f(a-\frac{h}{2})+ 23f(a+\frac{h}{2}) 
     + 24 f(a+\frac{3h}{2}) + \dots }{24} \\ 
     \frac{\dots + 24f(b-\frac{3h}{2}) + 23f(b-\frac{h}{2})
     + f(b+\frac{h}{2})}{24}
\end{multline}where $h\triangleq \frac{b-a}{N}$ is the step of integration,
$N$ is the number of points in which $f$ was evaluated, and $a$ and $b$ are the
integration limits (see \cref{fig:3rd_order_composite}).

\Cref{eqn:3rd_order_composite} can be rewritten as
\begin{equation}
    \int_a^b f(x) \,\dif x \approx M_{N-2}(f) + \Delta_3(f)
    \label{rule:3rd_order}
\end{equation}
where 
\begin{equation}
    M_{N-2}(f)\triangleq h\sum_{i=0}^{N-3}f(a+(i+\tfrac{1}{2})h)
\end{equation}
is the midpoint rule using $N-2$ points in the interval $[a,b]$, and
\begin{equation}
    \label{rule:3rd_order_correction}
    \Delta_3(f) \triangleq h \frac{f(a-h/2)-f(a+h/2) + f(b+h/2)-f(b-h/2)}{24}
\end{equation}
can be seen as a correction term to the midpoint rule.

\subsection{Derivation of an arbitrary order rule}
\label{sec:arbitrary_order}
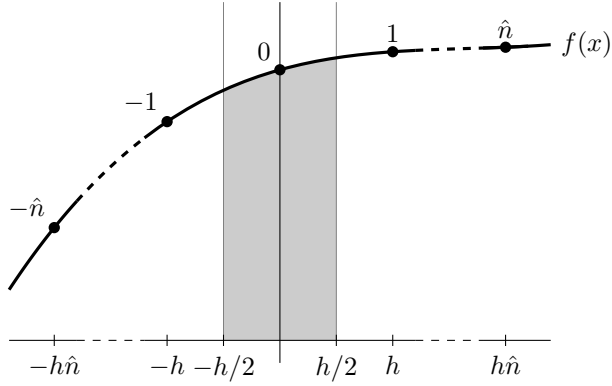
\begin{figure}
    \centering
    \begin{tikzpicture}
    [scale = 1.5,
     axis/.style = {thin},
     declare function  = {f(\x,\xa,\ya,\yo,\xb,\yb,\xd,\yd) = %
            ((\x-\xb)*\x*(\x-\xd))/((\xa-\xb)*\xa*(\xa-\xd))*\ya %
            + ((\x-\xb)*(\x-\xa)*(\x-\xd))/((-\xb)*(-\xa)*(-\xd))*\yo %
            + (\x*(\x-\xa)*(\x-\xd))/(\xb*(\xb-\xa)*(\xb-\xd))*\yb %
            + ((\x-\xb)*\x*(\x-\xa))/((\xd-\xb)*(\xd)*(\xd-\xa))*\yd %
            ;},
    ]
    \pgfmathsetmacro{\xa}{-2}
    \pgfmathsetmacro{\ya}{1}  %y(xa)
    \pgfmathsetmacro{\yo}{2.4} %y(0)
    \pgfmathsetmacro{\xb}{2}  
    \pgfmathsetmacro{\yb}{2.6} %y(xb)
    \pgfmathsetmacro{\xd}{1} 
    \pgfmathsetmacro{\yd}{2.56} %f(xd)
    \pgfmathsetmacro{\xmin}{-2.4}
    \pgfmathsetmacro{\xmax}{2.4}
    \pgfmathsetmacro{\ymin}{-.2}
    \pgfmathsetmacro{\ymax}{3}

    %f(x): 3rd degree polynomial
    \draw[domain=\xmin:-1.8][samples=20][very thick]
        plot (\x,{f(\x,\xa,\ya,\yo,\xb,\yb,\xd,\yd)}); 
    \draw[domain=-1.8:-1.2][samples=20][very thick, dashed]
        plot (\x,{f(\x,\xa,\ya,\yo,\xb,\yb,\xd,\yd)}); 
    \draw[domain=-1.2:1.2][samples=20][very thick]
        plot (\x,{f(\x,\xa,\ya,\yo,\xb,\yb,\xd,\yd)}); 
    \draw[domain=1.2:1.8][samples=20][very thick, dashed]
        plot (\x,{f(\x,\xa,\ya,\yo,\xb,\yb,\xd,\yd)}); 
    \draw[domain=1.8:2.2][samples=20][very thick]
        plot (\x,{f(\x,\xa,\ya,\yo,\xb,\yb,\xd,\yd)}); 
    \draw[domain=1.8:\xmax][samples=20][very thick]
        plot (\x,{f(\x,\xa,\ya,\yo,\xb,\yb,\xd,\yd)}) 
	[anchor=west] node {$f(x)$};
    
    %Interpolation points
    \foreach \p/\labeltext/\anchor in {(\xa,\ya)/-\hat{n}/south east,%
                    (-1,{f(-1,\xa,\ya,\yo,\xb,\yb,\xd,\yd)})/-1/south east,%
                    (0,\yo)/0/south east,%
                    (1,{f(1,\xa,\ya,\yo,\xb,\yb,\xd,\yd)})/1/south,%
                    (\xb,\yb)/\hat{n}/south%
                     }
    \fill \p circle [radius=.05] [anchor={\anchor}] node {$\labeltext$};
    
    \begin{scope}[on background layer]
        %Area corresponding to integration
        \fill[black!20!white][domain=-.5:.5]
             plot (\x,{f(\x,\xa,\ya,\yo,\xb,\yb,\xd,\yd)}) -- (.5,0) -- (-.5,0) --cycle;

	%Borders from the interval of integration
        \draw[help lines]
	    (-.5,0) -- (-.5,\ymax) (.5,0) -- (.5,\ymax);
    \end{scope}

    %Axis
	\draw[axis] (\xmin,0) -- (-1.8,0); %x axis
	\draw[axis, dashed] (-1.8,0)  -- (-1.2,0); %x axis
	\draw[axis] (-1.2,0)  -- (1.2,0); %x axis
	\draw[axis, dashed] (1.2,0)   -- (1.8,0); %x axis
	\draw[axis] (1.8,0) -- (\xmax,0); %x axis
	\draw
        [axis] (-2,.05) -- (-2,-.05)
        [below] node {$-h\hat{n}$}
        [axis] (-1,.05) -- (-1,-.05)
		[below] node {$-h$}
		[axis] (-.5,.05) -- (-.5,-.05)
		[below] node{$-h/2$}
		[axis] (.5,.05) -- (.5,-.05)
		[below] node{$h/2$}
		[axis] (1,.05) -- (1,-.05)
		[below] node{$h$}
        [axis] (2,.05) -- (2,-.05)
        [below] node {$h\hat{n}$};
	\draw[axis] (0, \ymin) -- (0, \ymax); %y axis

\end{tikzpicture}
    \caption{Quadrature rule of order $n=2\hat{n}+1$. The rule uses function
    evaluations on points $-h\hat{n}, \dots, h\hat{n}$, in order to approximate
    the integral in $[-h/2,h/2]$.}
    \label{fig:arbitrary_order}
\end{figure}

Extending the reasoning of the previous section, we can evaluate $f(x)$ in a
central point and $n-1$ symmetrically distributed side points, $f(x_k),\ x_k=kh,\
k=-\hat{n}, \dots, \hat{n}$, where $\hat{n} \triangleq \frac{n-1}{2}$
(see \cref{fig:arbitrary_order}). We can then construct the following interpolating
polynomial of degree $n-1$
\begin{equation}
    p_{n-1}(x) = \sum_{k=-\hat{n}}^{\hat{n}} \phi^n_k(x) f(x_k)
\end{equation}
where, by Lagrange's formula,
\begin{equation}
    \label{eqn:phi^n_k}
    \phi^n_k(x) = \prod_{\substack{i=-\hat{n} \\ i \neq k}}^{\hat{n}} 
    \frac{x - x_i}{x_k - x_i}
    = \prod_{\substack{i=-\hat{n} \\ i \neq k}}^{\hat{n}} 
    \frac{x-ih}{(k-i)h}
\end{equation}

Taking $p_{n-1}(x)$ as an approximation of $f(x)$ over $[-h/2,h/2]$, we can
evaluate the integral
\begin{equation}
    \int_{-h/2}^{h/2} f(x)\dif x \approx \int_{-h/2}^{h/2} p_{n-1} (x) \dif x 
\end{equation}
Denoting this result as the quadrature rule $L_n$, we have
\begin{equation}
    L_n\left(f(x)\right)=    h\sum_{k=-\hat{n}}^{\hat{n}} w^n_k f(x_k)
\end{equation}
where the normalized weights $w^n_k$, considering~\cref{eqn:phi^n_k}, are given by
\begin{equation}
    w^n_k \triangleq \frac{1}{h}\int_{-h/2}^{h/2} \phi^n_k(x)\dif x 
         = \frac{1}{h}\int_{-h/2}^{h/2} \prod_{\substack{i=-\hat{n} \\ i \neq k}}
    ^{\hat{n}}
    \frac{\frac{x}{h}-i}{k-i} \dif x
\end{equation}

In terms of the normalized variable $u \triangleq \tfrac{x}{h}$, we can write
\begin{equation}
    \label{eqn:weight}
    w^n_k = \int_{-1/2}^{1/2} \prod_{\substack{i=-\hat{n} \\ i \neq k}}^{\hat{n}}
    \frac{u-i}{k-i} \dif u
\end{equation}

For the $n$ point rule derived above, the related $N$ point composite rule
is
\begin{equation}
    \int_a^b f(x) \dif x \approx M_{N-(n-1)}(f) +\Delta_n(f)
\end{equation}
where
\begin{align}
    \Delta_n &= \left(\Delta_n\right)_{\text{outside}} -
    \left(\Delta_n\right)_{\text{inside}}\\
    \left(\Delta_n\right)_{\text{outside}} &= h \sum_{i=1}^{\hat{n}}
    \left(\sum_{k=i}^{\hat{n}} w^n_k\right)
    \left[f\left(a+(-i+\tfrac{1}{2})h\right) + f\left(b+(i-\tfrac{1}{2})h\right)\right]\\
    \left(\Delta_n\right)_{\text{inside}} &= h \sum_{i=1}^{\hat{n}}
    \left(\sum_{k=i}^{\hat{n}} w^n_k\right)
    \left[f\left(a+(i-\tfrac{1}{2})h\right) + f\left(b+(-i+\tfrac{1}{2})h\right)\right]
\end{align}

This time, the correction term $\Delta_n$ uses $\hat{n}$ function evaluations
outside the integration interval and $\hat{n}$ points within the interval. The
weighted function evaluations outside the interval are added to the midpoint
rule and the weighted function evaluations inside the interval are subtracted.

\subsubsection{Fifth order rule}
A particular  result can be obtained from the previous section by setting $n=5$:
\begin{equation}
    \int_{-h/2}^{h/2} f(x)\dif x \approx 
    h\frac{-17f(-2h) + 308f(-h) + 5178f(0) + 308f(h) - 17f(2h)}{5760}
\end{equation}

This simple rule can then be composed to yield
\begin{equation}
    \int_a^b f(x) \dif x \approx M_{N-4}(f) + \Delta_5(f)
\end{equation}
where $N$ is the number of points considered by the rule, $M_{N-4}(f)$ is
again the midpoint rule, and the correction term takes the form
\begin{multline*}
    \Delta_5 = h\left\{ 
    \frac{97}{1920} \left[f(a-\textstyle\frac{h}{2}) - f(a+\textstyle\frac{h}{2}) 
    - f(b-\textstyle\frac{h}{2}) + f(b+\textstyle\frac{h}{2})
    \right] 
    -\right. \\
    \left. 
    - \frac{17}{5760}\left[f(a-\textstyle\frac{3h}{2}) - f(a+\textstyle\frac{3h}{2})
    - f(b-\textstyle\frac{3h}{2}) + f(b+\textstyle\frac{3h}{2})
    \right] 
    \right\}
\end{multline*}

\section{Integration error}
\label{sec:integration_error}
\subsection{Notation}
For this section, a brief reminder of the used notation is in order. In this
paper, $n$ stands for the number of points used in the simple quadrature rule,
which, as seen in \cref{sec:arbitrary_order} is always odd. Due to the symmetry
of this rule family, it is convenient to define $\hat{n}=\frac{n-1}{2}$.
Finally, $N$ stands for the total number of points used by a composite rule.

Furthermore, we will denote a quadrature rule by the letter $L$ and the
exact definite integral by the cursive letter $\mathscr{L}$.
A superscript will be used to indicate the
number of points of the simple rule from which $L$ was derived, and a subscript
will indicate the total number of points used by $L$, since it may be composite.
For example, the notation $L^n_N$ stands for the $N$ point composite rule
derived from the $n$ point simple rule. The integration error is defined as
$R^n_N\triangleq L^n_N-\mathscr{L}$. In order to simplify
the notation, for simple rules we will drop the superscript so that $L_n
\triangleq L^n_n$.

\subsection{Error formula derivation}
\label{sec:error_derivation}
We begin by reminding ourselves that $L_n$ is an interpolatory rule, so it is
exact for polynomials of degree $n-1$ and bellow. In fact, we will prove that it
is exact for polynomials of degree $n$ as well, by introducing \cref{def:symmetric}
and then proceeding to prove \cref{lemma:exact}.

\begin{definition}
\label{def:symmetric}
We say that an integration rule is symmetric if and only if, for all weights
\[
    w_k = w_{-k}, \quad k = -\hat{n},\ldots,\hat{n}
\]
and for all function evaluation points
\[
     x_k = -x_{-k}, \quad\ k = -\hat{n},\ldots,\hat{n}
\]
\end{definition}

For the proposed quadrature rules, it is easy to see that $x_k = -x_{-k}$, since
$x_i \triangleq kh$. Moreover, from equation~\cref{eqn:weight}, substituting the
dummy variables $u$ and~$i$ for $-u$ and~$-i$ respectively, and using
multiplication's commutativity, we have 
\begin{multline}
    w_k = \int_{-1/2}^{1/2} \prod_{\substack{i=-\hat{n} \\ i \neq k}}^{\hat{n}}
    \frac{u-i}{k-i} \dif u =\\
    = \int_{-1/2}^{1/2} \prod_{\substack{i=-\hat{n} \\ i \neq -k}}^{\hat{n}}
    \frac{(-u)+i}{k+i} \dif (-u) =\\
    = \int_{-1/2}^{1/2} \prod_{\substack{i=-\hat{n} \\ i \neq -k}}^{\hat{n}}
    \frac{u-i}{-k-i} \dif u = w_{-k}
\end{multline}
Thus the proposed quadrature rules are symmetric

We can also prove the following result for symmetric interpolatory quadrature
rules

\begin{lemma}
\label{lemma:exact}
    A $n$ point symmetric interpolatory rule $L_n$ is exact for polynomials of up to
    degree $n$, i.e.\
\[
    R_{n}(p) = 0 \quad\forall\quad p \in \Pi_n,\Pi_{n-1},\dots,\Pi_0
\]
where $\Pi_n$ is the set of polynomials of degree $n$.

\begin{proof}
    Since $L_n$ is interpolatory and uses $n$ points, it is, by construction,
exact for all polynomials of degree $n-1$ or less. We will now prove that
this is also the case for polynomials of degree $n$. As usual, we define
$\hat{n}\triangleq (n-1)/2$.

    Let $p(x)= \sum_{i=0}^{\hat{n}}a_{2i}x^{2i} + \sum_{i=0}^{\hat{n}}
    a_{2i+1}x^{2i+1}$. Then
\[
    R_{2\hat{n}+1}(p(x)) =  R_{2\hat{n}+1}\left(\sum_{i=0}^{\hat{n}} a_{2i}x^{2i} \right) 
    + R_{2\hat{n}+1}\left(\sum_{i=0}^{\hat{n}} a_{2i+1}x^{2i+1}\right)
\]
The first term argument is a polynomial of degree $2\hat{n}=n-1$ and the rule is
interpolatory, so $R_{2\hat{n}+1}\left(\sum_{i=0}^{\hat{n}} a_{2i}x^{2i} \right) = 0$. 
    
The second term can be expanded as
\[
R_{2\hat{n}+1}(p(x)) = L_{2\hat{n}+1}\left(\sum_{i=0}^{\hat{n}} a_{2i+1}x^{2i+1}\right) - \mathscr{L}_{2\hat{n}+1}\left(\sum_{i=0}^{\hat{n}} a_{2i+1}x^{2i+1}\right)
\]
Since $\sum_{i=0}^{\hat{n}} a_{2i+1}x^{2i+1}$ is odd and $L$ is symmetric, both
terms in this equation are zero.
\end{proof}
\end{lemma}

Since our quadrature rules are symmetric and interpolatory, \cref{lemma:exact}
applies. This result allows the integration error to be derived by means of the
Peano's theory. 
Following the procedure suggested by David Ferguson \cite[example (b)]{davidferguson}, 
 the integration error is
\begin{equation}
    \label{eqn:error_peano_theory}
	R_{n}(f(x)) = f^{(n+1)}(\xi) R_{n}(\sfrac{x^{n+1}}{(n+1)!}) \quad \xi \in [-h,h]
\end{equation}
provided that $f$ is at least $n+1$ times differentiable.

A simple variable substitution allows $R_{n}$ to be computed in any arbitrary interval $[\alpha, \beta]$.
Defining  $u \triangleq \frac{x-a}{b-a}{\scriptstyle (\beta-\alpha)+\alpha}$, we have $\dif x=\frac{(b-a)}{\beta-\alpha} d\hat{x}$ and 

{
    \newcommand{\xu}{%
        \left({\scriptstyle \frac{u-\alpha}{\beta-\alpha}(b-a)+a}\right)%
    }
\begin{align}
    R_n(x^{n+1}) &= \int_a^b x^{n+1} \dif x - (b-a)\sum_{i} w_i x_i^{n+1} = \\
                 &= \int_\alpha^\beta \xu^{n+1} \dif u - (b-a) \sum_i w_i \xu^{n+1}
                    \nonumber
\end{align}
}
where $w_i$ are the $n$ normalized weights.

Since the rule is exact for polynomials up to degree n, we can drop those when expanding the terms in $u$
\begin{equation}
    R_n(x^{n+1}) = (b-a)^{n+2}
    \left[
        \frac{\int_\alpha^\beta u^{n+1} \dif u - (\beta-\alpha)\sum_i w_i u_i}{(\beta-\alpha)^{n+2}}
    \right]
    \quad \forall \, \alpha,\beta
\end{equation}

This allows us to define the normalized integration error
\begin{equation}
    \label{eqn:normalized_error}
    \hat{R}_n \triangleq \frac{R_n(\sfrac{x^{n+1}}{(n+1)!})}{h^{n+2}}
\end{equation}
which is constant for a given rule. This result is valid in general for any rule with error given by \cref{eqn:error_peano_theory} with fixed number of steps per integration interval.%

The integration error can now be written as
\begin{equation}
    R_{n}(f(x)) = \hat{R}_n h^{n+2} f^{(n+1)}(\xi) \quad \xi \in [-h,h]
\end{equation}

This proves that the $n$ point simple quadrature rule is indeed of $n$th
order, that is, it integrates polynomials of degree up to $n$ exactly.

\subsection{Integration error of the composite rule}

The composite rule is derived from applying the simple rule to $M$
sub-intervals, one for each step. Therefore its error is given by
\begin{equation}
\label{eqn:composite_error1}
R_N=\sum_{i=1}^{M} \hat{R}_n h^{n+2} f^{(n+1)}(\xi_i)
\end{equation}
where $\xi_i \in [a+(i-\frac{3}{2})h,a+(i+\frac{1}{2})h]$, i.e.\ $\xi$ is in
a radius $h$ ball centered in each integration point. The integration step is
given by $h=\frac{b-a}{M}$. We can rewrite \cref{eqn:composite_error1} as
\begin{equation}
    \label{eqn:composite_error2}
    R_N=\hat{R}_n h^{n+2} M\sum_{i=1}^{M}\frac{f^{(n+1)}(\xi_i)}{M}
\end{equation}
If we assume $f^{(n+1)}$ continuous, there is a $\xi \in [a-h/2, b+h/2]$ that
satisfies 
\begin{equation}
f^{(n+1)}(\xi)=\sum_{i=1}^{M} \frac{f^{(n+1)}(\xi_i)}{M}
\end{equation}
so \cref{eqn:composite_error2} becomes
\begin{equation}
R_N=\hat{R}_n M h^{n+2} f^{(n+1)}(\xi)
\end{equation}
Substituting $h$,
\begin{equation}
R_N = \hat{R}_n(b-a)h^{n+1}f^{(n+1)}(\xi)
\end{equation}
or,
\begin{equation}
    \label{eqn:composite_rule3}
R_N = \hat{R}_n \frac{(b-a)^{n+2} f^{(n+1)}(\xi)}{M^{n+1}}
\end{equation}

Note that the simple $n$th order rule uses one point inside each integration step and
$n-1$ points outside. Since the points outside the interval are also used by one or
more neighboring steps, the composite rule uses $n-1$ points outside the full
integration interval, namely $(n-1)/2$ points before the beginning 
and $(n-1)/2$ points after the end. In addition to that, one point for
each step is used. Therefore the total number of points is given by $N=M+n-1$.
Solving for $M$, we have
\begin{equation}
    M=N-(n-1)
\end{equation}
Substituting this result in \cref{eqn:composite_rule3}, we have our final error
formula,
\begin{equation}
   R_N = \hat{R}_n \frac{(b-a)^{n+2} f^{(n+1)}(\xi)}{[N-(n-1)]^{n+1}} \quad \xi \in [a-(n-2)h/2, b+(n-2)h/2]
\end{equation}

\section{Modified rules}
\label{sec:modified_rules}
\subsection{Formulae using only points in the integration interval}
\label{sec:interval_rule}

The major disadvantage of the family of rules presented so far in relation to
the traditional Newton-Cotes family is that the former uses points
outside the integration interval: in some cases these points may simply not be available for
evaluation. Nevertheless, a simple adaptation of the end-steps can avoid this
issue entirely. The modification will be derived here only for the third order
rule for brevity, but the reader will find it easy to apply the same reasoning
to any other case.

\begin{figure}
\centering
\begin{tikzpicture}
    [scale = 3,
        axis/.style = {thin}
    ]
    \pgfmathsetmacro{\ya}{.5}  %y(-1)
    \pgfmathsetmacro{\yo}{1.2} %y(0)
    \pgfmathsetmacro{\yb}{1.3} %y(1)
    \pgfmathsetmacro{\yd}{1.28} %p_2(.5)
    \pgfmathsetmacro{\xmin}{-1.2}
    \pgfmathsetmacro{\xmax}{1.2}
    \pgfmathsetmacro{\ymin}{-.2}
    \pgfmathsetmacro{\ymax}{1.5}

    \pgfmathsetmacro{\yA}{0.97}  %y(a)

    %f(x): 3rd degree polynomial
    \draw[domain=\xmin:\xmax][samples=100][very thick]
        plot (\x,{-(\x-1)*\x*(\x-.5)/3*\ya+ 2*(\x-1)*(\x+1)*(\x-.5)*\yo + \x*(\x+1)*(\x-.5)*\yb + -(\x-1)*\x*(\x+1)*4/1.5*\yd}) 
	[anchor=south west] node {$f(x)$};
    
    %Interpolating parabola
    \draw[domain=\xmin:\xmax][samples=100][semithick]
        plot (\x,{(\x-1)*\x/0.75*\yA -(\x-1)*(\x+0.5)/0.5*\yo + \x*(\x+.5)/1.5*\yb}) 
	[anchor=north west] node {$p_2(x)$};

    %Interpolation points
    \foreach \p in {(-0.5,\yA), (0,\yo), (1,\yb)}
    \fill \p circle [radius=.025];
    
    \begin{scope}[on background layer]
        %Area corresponding to integration
        \fill[black!20!white][domain=-.5:.5]
             plot (\x,{(\x-1)*\x/0.75*\yA -(\x-1)*(\x+0.5)/0.5*\yo + \x*(\x+.5)/1.5*\yb}) -- (.5,0) -- (-.5,0) --cycle;

	%Borders from the interval of integration
        \draw[help lines]
	    (-.5,0) -- (-.5,\ymax) (.5,0) -- (.5,\ymax);
    \end{scope}

    %Axis
	\draw[axis] (\xmin,0) -- (\xmax,0); %x axis
    \draw [axis,dashed] (\xmax,0) -- (\xmax+0.35,0);
    %\foreach \x in {-1, -1/2, 1/2, 1}
    %    \draw[axis] (\x,.05) -- (\x,-.05) %x axis ticks
	%    [below] node {$\x$}; %x axis labels
	\draw[axis] (-1,.05) -- (-1,-.05)
		[axis] (-.5,.05) -- (-.5,-.05)
		[below] node{$a$}
		[axis] (0,.05) -- (0,-.05)
        [below] node{$a+\frac{1}{2}h$}
		[axis] (.5,.05) -- (.5,-.05)
		[below] node{$a+h$}
		[axis] (1,.05) -- (1,-.05)
        [below] node{$a+\frac{3}{2}h$};

\end{tikzpicture}
\caption{Modified first step for third order quadrature. The interpolating
    polynomial $p_2$, which is used to approximate $f$ is constructed from
    points $a,a+\frac{1}{2}h,a+\frac{3}{2}h$ in order to avoid points outside
    the interval. The integration is carried from $a$ to $a+h$ as usual. Compare
    with \cref{fig:3rd_order_simple}.}
\label{fig:interval_rule}
\end{figure}
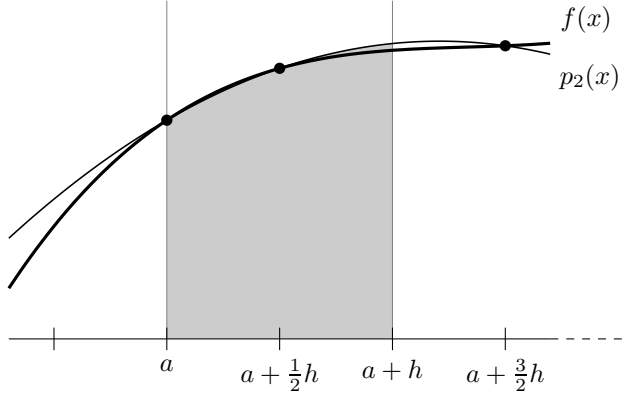

This modification consists on modifying the interpolating polynomial of the
end-steps in such a fashion to evaluate the end-point instead of the point
outside the integration interval, as shown in \cref{fig:interval_rule} for the
first step.

In order to simplify the derivation for the first step, we define the auxiliary
variable 
\begin{equation} 
    u \triangleq \frac{x - (a+\frac{1}{2}h)}{h}
\end{equation}
Doing so is equivalent to assuming, without loss of generality, that
 $a+\frac{1}{2}h$ is the origin and that the integration step is unitary.
The interpolating polynomial $p_2(u)$ is 
\begin{equation} 
    p_2(u) =
    \frac{u(u-1)}{(-\tfrac{1}{2})(-\tfrac{3}{2})}f(a)
    + \frac{(u+\tfrac{1}{2})(u-1)}{(\tfrac{1}{2})(-1)} f(a+\tfrac{1}{2}h)
    + \frac{(u+\tfrac{1}{2})u}{(\tfrac{3}{2})(1)} f(a+\tfrac{3}{2}h)
\end{equation}
Integrating, we have
\begin{multline}
    \label{eqn:interval_first}
    \int_a^{a+h}f(x) \dif x 
    \approx \int_a^{a+h}p_2(x) \dif x = \\
    = \int_{-\tfrac{1}{2}}^{\tfrac{1}{2}} p_2(u) h \dif u
    = \frac{2f(a)+15f(a+\tfrac{1}{2}h)+f(a+\tfrac{3}{2}h)}{18}h
\end{multline}

The result for the last step is symmetrical to \ref{eqn:interval_first}, i.e.\
\begin{equation}
    \label{eqn:interval_last}
    \int_{b-h}^b f(x)\dif x 
    \approx \frac{2f(b)+15f(b-\tfrac{1}{2}h)+f(b-\tfrac{3}{2}h)}{18}h
\end{equation}

Finally, for the remaining steps, the usual (\ref{sec:3rd_order}) rule is used,
i.e.\
\begin{multline}
    \label{eqn:interval_middle}
    \int_{a+h}^{b-h} f(x)\dif x = h\sum_{i=0}^{N-5}f(a+(i+\tfrac{3}{2})h) +\\
    + h\frac{f(a+\tfrac{1}{2}h)-f(a+\tfrac{3}{2}h)
    + f(b-\tfrac{1}{2}h)-f(a-\tfrac{3}{2}h)}{12}
\end{multline}

Summing \cref{eqn:interval_first,eqn:interval_last,eqn:interval_middle},
we get the final rule
\begin{multline}
    \label{eqn:interval}
    \int_a^b f(x)\dif x \approx M_{N-2} +\\+ h\frac{f(a)+f(b)}{9} 
    - h\frac{3f(a+\tfrac{1}{2}h)+f(a+\tfrac{3}{2}h)+f(b-\tfrac{3}{2}h)+f(b-\tfrac{1}{2}h)}{36}
\end{multline}

This modified rule is still of third order, since it is symmetric and
interpolatory in a stepwise sense. The proof for this more
general\footnote{In this case, each step of the composite rule is interpolatory.
Therefore, when calculating the error in each step for a polynomial one degree above that of
the interpolating polynomial, which is even by construction, the even part of the error
will vanish because it is of the order of the interpolation, while the odd part
will vanish because the quadrature rule is symmetric.} case is
entirely analogous to that made in \cref{lemma:exact}.

\subsection{Formulae using endpoint derivatives}
\label{sec:derivative_rule}

Another related integration formula can be derived by a similar reasoning, but
using end-point derivatives instead of points outside the integration interval.
This formula was derived by \cite{derivative_rule}, albeit with a different reasoning.

As a motivation, consider the third order correction term
\[
    \Delta_3=\frac{h^2}{24}
    \left[
        \frac{f(b+h/2)-f(b-h/2)}{h} - \frac{f(a+h/2)-f(a-h/2)}{h} 
    \right]
\]
Assuming a small $h$, we have 
\begin{equation}
    \label{eqn:delta3prime}
    \Delta_3 \approx \Delta_3' = \frac{h^2}{24} (f'(b)-f'(a))
\end{equation}

Now we will show that the correction term $\Delta_3'$ 
summed to the midpoint rule yields a third order integration rule. 
Assume $f:[-h/2,h/2] \rightarrow \mathbb{R}$ differentiable with a continuous derivative. 
Approximate it by a quadratic polynomial $p_2(x)$ such that 
$p_2'(-h/2)=f'(-h/2)$, $p_2'(h/2)=f'(h/2)$ and $p_2(0)=f(0)$, i.e.\
\begin{equation} 
    p_2(x) = \frac{x^2-hx}{-2h}f'(-h/2)+f(0)+\frac{x^2+hx}{2h}f'(h/2)
\end{equation}

Integrating $p_2(x)$, we get
\begin{equation}
    \int_{-h/2}^{h/2}f(x)\,dx \approx
    \int_{-h/2}^{h/2} p_2(x)\,dx = \frac{h^2}{24}(f'(h/2)-f'(-h/2)) + h f(0)
\end{equation}
This rule can be readily composed to yield
\begin{equation}
    \int_{a}^{b}f(x)\,dx \approx M_{n-2} + \Delta_3'
\end{equation}
where
\begin{equation}
   \Delta_3'=\frac{h^2}{24}(f'(b)-f'(a))
\end{equation}
Notice that the derivatives taken inside the integration interval cancel out nicely.

One interesting property of the correction term $\Delta_3'$ is that it can be used to estimate the integration step for a given error \emph{a priori}. If we take the correction term to be an approximation of the error, denoted
$\bar{R}$, the integration step for a desired amount of error should be
\begin{equation} h=\sqrt{\left|\frac{24 \bar{R}}{f'(b)-f'(a)}\right|}
\end{equation} as long as $f'(b) \neq f'(a)$.

\section{Results and discussion}
\label{sec:results}
In~\cref{sec:numerical_results} we show several practical results of the
application of the rules derived in a set of representative functions selected
mostly from \cite{casaletto1969comparison}.

In general, the third order and the derivative rule compare favorably to
Simpson's rule, which is used as a benchmark. The modified third order rule with
no points outside the integration interval and no use of derivatives, derived in
\cref{sec:interval_rule}, also has errors on par with Simpson's rule.

\begin{table}
    \centering
    \caption{Normalized weights for rules up to ninth order.} 
    \label{tbl:weights}
    {
        \renewcommand{\arraystretch}{1.3}
        {
\newcolumntype{C}{>{$}c<{$}}
\begin{tabular}{CCCCCC}
\toprule 
2n+1 	& \multicolumn{5}{l}{$\qquad w_i, \quad i=0,\dots,n$} \\ 
\midrule 
1 &	 \frac{1}{1}\\
3 &	 \frac{11}{12}	&\frac{1}{24}\\
5 &	 \frac{863}{960}	&\frac{77}{1440}	&\frac{-17}{5760}\\
7 &	 \frac{215641}{241920}	&\frac{6361}{107520}	&\frac{-281}{53760}	&\frac{367}{967680}\\
9 &	 \frac{41208059}{46448640}	&\frac{3629953}{58060800}	&\frac{-801973}{116121600}	&\frac{49879}{58060800}	&\frac{-27859}{464486400}\\
\bottomrule
\end{tabular}
}

    }
\end{table}

\begin{figure}
    \centering
    \textbf{Numerical stability}
    \input{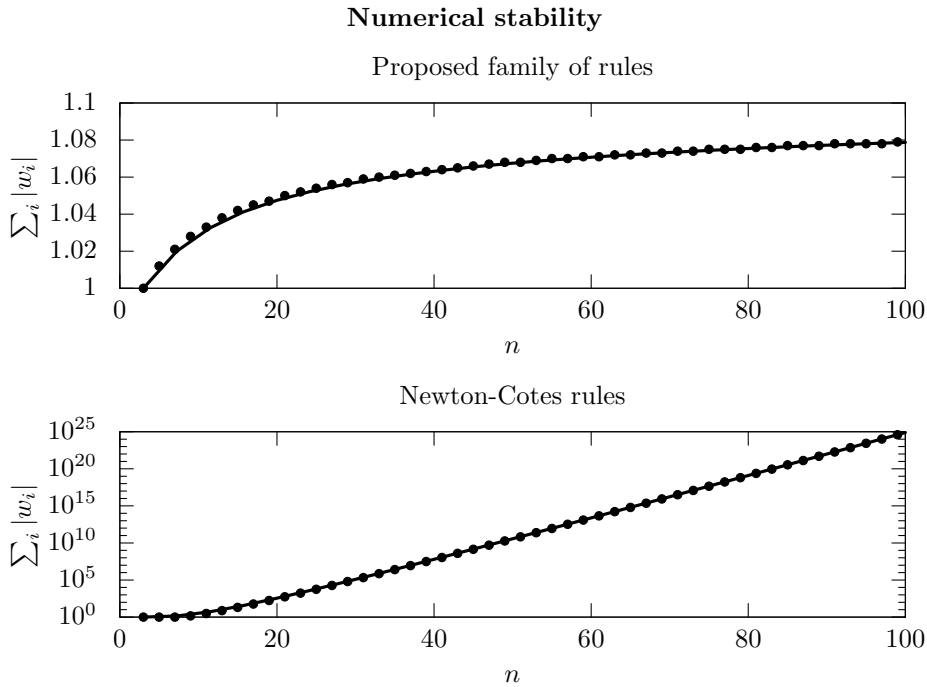}
    \label{plot:stability}
    \caption
    {
    Sum of normalized absolute rule weights for the proposed family and Newton-Cotes rules. The weights were normalized so that $\sum_i w_i=1$. For the
    proposed family, the sum keeps quite close to unity, while it diverges
    quickly for Newton-Cotes formulae.
    }
\end{figure}

\Cref{tbl:weights} shows the tabulated rule weights for the proposed rules up to
ninth order. The presence of negative weights starting from the fifth order rule
suggests that these rules might become numerically unstable as you increase the
order, but that is not the case. \Cref{plot:stability} shows that the sum of
weights apparently has asymptotic behavior with an asymptote smaller than
1.1. Further testing has shown that this behavior holds at least up to \nth{420}
order, while the sum of weights for Newton-Cotes formulae diverges
exponentially.

The general results obtained for integration error in
\cref{sec:integration_error} can be particularized. For the third and fifth
order rule we get the following error formulae for the simple rules
\begin{align}
    R_3(f(x)) &= \frac{-17}{5760}\,h^{n+2} f^{(n+1)}(\xi)   & \xi \in [-h,h]
\\
    R_5(f(x)) &= \frac{367}{967680}\,h^{n+2} f^{(n+1)}(\xi) & \xi \in [-h,h]
\end{align}
and for the composite rules
\begin{align}
    R_N^3 &= -\frac{17}{5760}\frac{(b-a)^5 f^{(4)}(\xi)}{(N-2)^4}   & \xi \in [a-h/2, b+h/2] \\
    R_N^5 &= \frac{367}{967680}\frac{(b-a)^7 f^{(6)}(\xi)}{(N-4)^6} & \xi \in [a-3h/2, b+3h/2]
\end{align}

For the derivative rule, \cref{sec:derivative_rule}, we have
\begin{equation}
    R_N = -\frac{7}{5760}\frac{(b-a)^5 f^{(4)}(\xi)}{(N-2)^4} \qquad \xi \in [a, b]
\end{equation}
counting the derivative evaluations as function evaluations. 
It is slightly better than the regular third order rule.

We also show values of the normalized integration error $\hat{R}_n$ for rules up to ninth rule in \cref{tbl:Kn}

\begin{table}
\caption{Normalized integration errors for the proposed rules up to ninth order} 
{
\renewcommand{\arraystretch}{1.3}
\label{tbl:Kn}
\[ \begin{array}{ccc}
\toprule
n	& {\hat{R}_n} {\text{\tiny (exact)}}	& {\hat{R}_n} {\text{\tiny (approx.)}} \\
\midrule
1	& \frac{1}{24}	& 0.0416667\\3	& \frac{-17}{5760}	& -0.00295139\\5	&
\frac{367}{967680}	& 0.000379258\\7	& \frac{-27859}{464486400}	&
-5.99781\cdot10^{-5}\\9	& \frac{1295803}{122624409600}	& 1.05673\cdot 10^{-5}\\\bottomrule
\end{array} \] 

}
\end{table}

\subsection{Error formula comparison with Newton-Cotes integration rules} 

In order to compare the family of rule presented herein with Newton-Cotes
rules in an analytical manner, it is interesting to consider both error 
formulae in a common framework. For that we define the global normalized integration error as 
\begin{gather}
    \label{eqn:global_normalized_error}
    \hat{R}_N^n \triangleq  \frac{R_N^n(\sfrac{x^{n+1}}{(n+1)!})}{(b-a)^{n+2}} \\
    \makebox[\textwidth][l]{so that}\nonumber \\ 
    R_N^n = \hat{R}_N^n (b-a)^{n+2} f^{(n+1)}(\xi)
\end{gather}
where $n$ and $M$ retain their meaning of rule order and number sub-intervals respectively, 
$a$ and $b$ are the integration bounds, and $\xi \in
[a-(n-2)h/2,b+(n-2)h/2]$ in the general case, or more strictly $\xi \in
[a,b]$ for Newton-Cotes rules. $\hat{R}_n$ is the normalized integration error,
 as defined in \cref{eqn:normalized_error}.

From \cref{eqn:composite_rule3},
\begin{equation}
    \label{eqn:global_normalized_error2}
    \hat{R}_N^n = \frac{\hat{R}_n}{M^{n+1}}
\end{equation}

We should also write the global normalized integration error as function of
the number of points evaluated rather than number of sub-intervals. This is
because the number of sub-intervals is $N-1$ for Newton-Cotes rules, 
while  it is $N-(n-1)$ for the family of rules presented here.
Substituting $M$ in \cref{eqn:global_normalized_error2} and dropping the super and subscripts, we have
\begin{align}
    \hat{R}_\text{new family}   &= \frac{(\hat{R}_n)_\text{new family}}{[N-(n-1)]^{n+1}}\\
    \hat{R}_\text{Newton-Cotes} &= \frac{(\hat{R}_n)_\text{Newton-Cotes}}{(N-1)^{n+1}}
\end{align}

Since $\hat{R}$ is constant for a given rule with known order and number of
points used, it can be used as a figure of merit to access the relative
errors between rules of the same order.

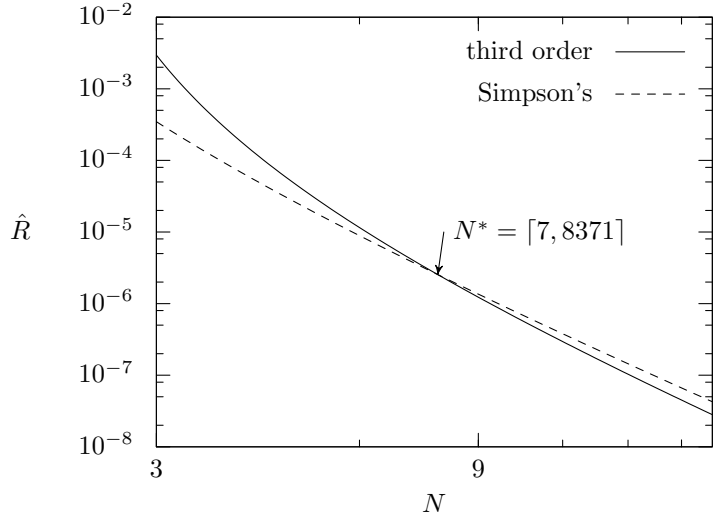
\begin{figure}
\label{plot:third_order_vs_simpson}
	\centering
	\begin{tikzpicture}[gnuplot]
%% generated with GNUPLOT 5.0p3 (Lua 5.1; terminal rev. 99, script rev. 100)
%% Sat 22 Sep 2018 07:24:32 AM -03
\gpmonochromelines
\path (0.000,0.000) rectangle (10.000,7.070);
\gpcolor{color=gp lt color border}
\gpsetlinetype{gp lt border}
\gpsetdashtype{gp dt solid}
\gpsetlinewidth{1.00}
\draw[gp path] (2.056,0.985)--(2.236,0.985);
\draw[gp path] (9.447,0.985)--(9.267,0.985);
\node[gp node right] at (1.872,0.985) {$10^{-8}$};
\draw[gp path] (2.056,1.272)--(2.146,1.272);
\draw[gp path] (9.447,1.272)--(9.357,1.272);
\draw[gp path] (2.056,1.651)--(2.146,1.651);
\draw[gp path] (9.447,1.651)--(9.357,1.651);
\draw[gp path] (2.056,1.845)--(2.146,1.845);
\draw[gp path] (9.447,1.845)--(9.357,1.845);
\draw[gp path] (2.056,1.938)--(2.236,1.938);
\draw[gp path] (9.447,1.938)--(9.267,1.938);
\node[gp node right] at (1.872,1.938) {$10^{-7}$};
\draw[gp path] (2.056,2.224)--(2.146,2.224);
\draw[gp path] (9.447,2.224)--(9.357,2.224);
\draw[gp path] (2.056,2.604)--(2.146,2.604);
\draw[gp path] (9.447,2.604)--(9.357,2.604);
\draw[gp path] (2.056,2.798)--(2.146,2.798);
\draw[gp path] (9.447,2.798)--(9.357,2.798);
\draw[gp path] (2.056,2.890)--(2.236,2.890);
\draw[gp path] (9.447,2.890)--(9.267,2.890);
\node[gp node right] at (1.872,2.890) {$10^{-6}$};
\draw[gp path] (2.056,3.177)--(2.146,3.177);
\draw[gp path] (9.447,3.177)--(9.357,3.177);
\draw[gp path] (2.056,3.556)--(2.146,3.556);
\draw[gp path] (9.447,3.556)--(9.357,3.556);
\draw[gp path] (2.056,3.751)--(2.146,3.751);
\draw[gp path] (9.447,3.751)--(9.357,3.751);
\draw[gp path] (2.056,3.843)--(2.236,3.843);
\draw[gp path] (9.447,3.843)--(9.267,3.843);
\node[gp node right] at (1.872,3.843) {$10^{-5}$};
\draw[gp path] (2.056,4.130)--(2.146,4.130);
\draw[gp path] (9.447,4.130)--(9.357,4.130);
\draw[gp path] (2.056,4.509)--(2.146,4.509);
\draw[gp path] (9.447,4.509)--(9.357,4.509);
\draw[gp path] (2.056,4.703)--(2.146,4.703);
\draw[gp path] (9.447,4.703)--(9.357,4.703);
\draw[gp path] (2.056,4.796)--(2.236,4.796);
\draw[gp path] (9.447,4.796)--(9.267,4.796);
\node[gp node right] at (1.872,4.796) {$10^{-4}$};
\draw[gp path] (2.056,5.082)--(2.146,5.082);
\draw[gp path] (9.447,5.082)--(9.357,5.082);
\draw[gp path] (2.056,5.462)--(2.146,5.462);
\draw[gp path] (9.447,5.462)--(9.357,5.462);
\draw[gp path] (2.056,5.656)--(2.146,5.656);
\draw[gp path] (9.447,5.656)--(9.357,5.656);
\draw[gp path] (2.056,5.748)--(2.236,5.748);
\draw[gp path] (9.447,5.748)--(9.267,5.748);
\node[gp node right] at (1.872,5.748) {$10^{-3}$};
\draw[gp path] (2.056,6.035)--(2.146,6.035);
\draw[gp path] (9.447,6.035)--(9.357,6.035);
\draw[gp path] (2.056,6.414)--(2.146,6.414);
\draw[gp path] (9.447,6.414)--(9.357,6.414);
\draw[gp path] (2.056,6.609)--(2.146,6.609);
\draw[gp path] (9.447,6.609)--(9.357,6.609);
\draw[gp path] (2.056,6.701)--(2.236,6.701);
\draw[gp path] (9.447,6.701)--(9.267,6.701);
\node[gp node right] at (1.872,6.701) {$10^{-2}$};
\draw[gp path] (2.056,0.985)--(2.056,1.165);
\draw[gp path] (2.056,6.701)--(2.056,6.521);
\node[gp node center] at (2.056,0.677) {$3$};
\draw[gp path] (4.756,0.985)--(4.756,1.075);
\draw[gp path] (4.756,6.701)--(4.756,6.611);
\draw[gp path] (6.336,0.985)--(6.336,1.075);
\draw[gp path] (6.336,6.701)--(6.336,6.611);
\draw[gp path] (7.457,0.985)--(7.457,1.075);
\draw[gp path] (7.457,6.701)--(7.457,6.611);
\draw[gp path] (8.326,0.985)--(8.326,1.075);
\draw[gp path] (8.326,6.701)--(8.326,6.611);
\draw[gp path] (9.037,0.985)--(9.037,1.075);
\draw[gp path] (9.037,6.701)--(9.037,6.611);
\draw[gp path] (6.336,0.985)--(6.336,1.165);
\draw[gp path] (6.336,6.701)--(6.336,6.521);
\node[gp node center] at (6.336,0.677) {$9$};
\draw[gp path] (9.037,0.985)--(9.037,1.075);
\draw[gp path] (9.037,6.701)--(9.037,6.611);
\draw[gp path] (2.056,6.701)--(2.056,0.985)--(9.447,0.985)--(9.447,6.701)--cycle;
\node[gp node left] at (5.877,3.843) {$N^*=\lceil 7,8371 \rceil$};
\draw[gp path,->](5.877,3.843)--(5.797,3.276);
\node[gp node center] at (0.246,3.843) {$\hat{R}$};
\node[gp node center] at (5.751,0.215) {$N$};
\node[gp node right] at (7.979,6.240) {third order};
\draw[gp path] (8.163,6.240)--(9.079,6.240);
\draw[gp path] (2.056,6.196)--(2.131,6.103)--(2.205,6.013)--(2.280,5.926)--(2.355,5.841)%
  --(2.429,5.760)--(2.504,5.681)--(2.579,5.603)--(2.653,5.528)--(2.728,5.455)--(2.803,5.384)%
  --(2.877,5.314)--(2.952,5.246)--(3.027,5.179)--(3.101,5.114)--(3.176,5.050)--(3.251,4.987)%
  --(3.325,4.925)--(3.400,4.865)--(3.474,4.805)--(3.549,4.746)--(3.624,4.689)--(3.698,4.632)%
  --(3.773,4.576)--(3.848,4.521)--(3.922,4.467)--(3.997,4.413)--(4.072,4.360)--(4.146,4.308)%
  --(4.221,4.256)--(4.296,4.205)--(4.370,4.154)--(4.445,4.104)--(4.520,4.055)--(4.594,4.006)%
  --(4.669,3.958)--(4.744,3.910)--(4.818,3.863)--(4.893,3.816)--(4.968,3.769)--(5.042,3.723)%
  --(5.117,3.677)--(5.192,3.632)--(5.266,3.587)--(5.341,3.542)--(5.416,3.498)--(5.490,3.454)%
  --(5.565,3.410)--(5.640,3.367)--(5.714,3.324)--(5.789,3.281)--(5.863,3.239)--(5.938,3.196)%
  --(6.013,3.154)--(6.087,3.113)--(6.162,3.071)--(6.237,3.030)--(6.311,2.989)--(6.386,2.948)%
  --(6.461,2.908)--(6.535,2.868)--(6.610,2.828)--(6.685,2.788)--(6.759,2.748)--(6.834,2.708)%
  --(6.909,2.669)--(6.983,2.630)--(7.058,2.591)--(7.133,2.552)--(7.207,2.514)--(7.282,2.475)%
  --(7.357,2.437)--(7.431,2.399)--(7.506,2.360)--(7.581,2.323)--(7.655,2.285)--(7.730,2.247)%
  --(7.805,2.210)--(7.879,2.172)--(7.954,2.135)--(8.029,2.098)--(8.103,2.061)--(8.178,2.024)%
  --(8.252,1.987)--(8.327,1.951)--(8.402,1.914)--(8.476,1.878)--(8.551,1.841)--(8.626,1.805)%
  --(8.700,1.769)--(8.775,1.733)--(8.850,1.697)--(8.924,1.661)--(8.999,1.626)--(9.074,1.590)%
  --(9.148,1.554)--(9.223,1.519)--(9.298,1.483)--(9.372,1.448)--(9.447,1.413);
\node[gp node right] at (7.979,5.678) {Simpson's};
\gpsetdashtype{gp dt 2}
\draw[gp path] (8.163,5.678)--(9.079,5.678);
\draw[gp path] (2.056,5.311)--(2.131,5.263)--(2.205,5.216)--(2.280,5.170)--(2.355,5.124)%
  --(2.429,5.078)--(2.504,5.033)--(2.579,4.988)--(2.653,4.943)--(2.728,4.899)--(2.803,4.855)%
  --(2.877,4.812)--(2.952,4.768)--(3.027,4.725)--(3.101,4.683)--(3.176,4.640)--(3.251,4.598)%
  --(3.325,4.556)--(3.400,4.515)--(3.474,4.473)--(3.549,4.432)--(3.624,4.391)--(3.698,4.350)%
  --(3.773,4.310)--(3.848,4.270)--(3.922,4.230)--(3.997,4.190)--(4.072,4.150)--(4.146,4.110)%
  --(4.221,4.071)--(4.296,4.032)--(4.370,3.993)--(4.445,3.954)--(4.520,3.916)--(4.594,3.877)%
  --(4.669,3.839)--(4.744,3.801)--(4.818,3.763)--(4.893,3.725)--(4.968,3.687)--(5.042,3.650)%
  --(5.117,3.612)--(5.192,3.575)--(5.266,3.538)--(5.341,3.501)--(5.416,3.464)--(5.490,3.427)%
  --(5.565,3.390)--(5.640,3.353)--(5.714,3.317)--(5.789,3.280)--(5.863,3.244)--(5.938,3.208)%
  --(6.013,3.172)--(6.087,3.136)--(6.162,3.100)--(6.237,3.064)--(6.311,3.028)--(6.386,2.993)%
  --(6.461,2.957)--(6.535,2.921)--(6.610,2.886)--(6.685,2.851)--(6.759,2.815)--(6.834,2.780)%
  --(6.909,2.745)--(6.983,2.710)--(7.058,2.675)--(7.133,2.640)--(7.207,2.605)--(7.282,2.571)%
  --(7.357,2.536)--(7.431,2.501)--(7.506,2.467)--(7.581,2.432)--(7.655,2.398)--(7.730,2.363)%
  --(7.805,2.329)--(7.879,2.295)--(7.954,2.260)--(8.029,2.226)--(8.103,2.192)--(8.178,2.158)%
  --(8.252,2.124)--(8.327,2.090)--(8.402,2.056)--(8.476,2.022)--(8.551,1.988)--(8.626,1.954)%
  --(8.700,1.921)--(8.775,1.887)--(8.850,1.853)--(8.924,1.819)--(8.999,1.786)--(9.074,1.752)%
  --(9.148,1.719)--(9.223,1.685)--(9.298,1.652)--(9.372,1.618)--(9.447,1.585);
\gpsetdashtype{gp dt solid}
\draw[gp path] (2.056,6.701)--(2.056,0.985)--(9.447,0.985)--(9.447,6.701)--cycle;
%% coordinates of the plot area
\gpdefrectangularnode{gp plot 1}{\pgfpoint{2.056cm}{0.985cm}}{\pgfpoint{9.447cm}{6.701cm}}
\end{tikzpicture}
%% gnuplot variables
    \caption{$\hat{R}$ constant comparison between the third order rule and Simpson's
    rule.}
\end{figure}

In particular, for the third order rule, $|\hat{R}|=\frac{17}{5760(N-2)^4}$, while,
for Simpson's rule, $|\hat{R}|=\frac{1}{180(N-1)^4}$. This values have been plotted in
\cref{plot:third_order_vs_simpson} as function of $N$. Simpson's rule
starts with a smaller error constant, but for $N \geq 8$ the third order
rule's error constant gets smaller. Asymptotically, the ratio $\hat{R}_\text{third
order}/\hat{R}_\text{Simpson}$ goes to $17/32 \approx 0.53$.

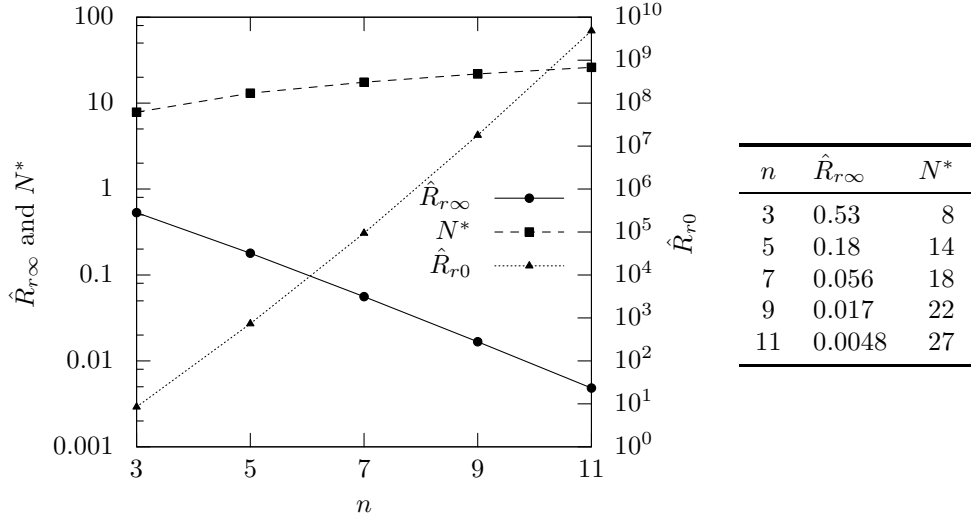
\begin{figure}
    \label{fig:utschVSnewtoncotes}
    \begin{minipage}{.755\textwidth}
      {\begin{tikzpicture}[gnuplot]
%% generated with GNUPLOT 5.0p3 (Lua 5.1; terminal rev. 99, script rev. 100)
%% Sat 22 Sep 2018 07:28:15 AM -03
\gpmonochromelines
\path (0.000,0.000) rectangle (9.700,7.070);
\gpcolor{color=gp lt color border}
\gpsetlinetype{gp lt border}
\gpsetdashtype{gp dt solid}
\gpsetlinewidth{1.00}
\draw[gp path] (1.688,0.985)--(1.868,0.985);
\node[gp node right] at (1.504,0.985) {$0.001$};
\draw[gp path] (1.688,1.329)--(1.778,1.329);
\draw[gp path] (1.688,1.784)--(1.778,1.784);
\draw[gp path] (1.688,2.017)--(1.778,2.017);
\draw[gp path] (1.688,2.128)--(1.868,2.128);
\node[gp node right] at (1.504,2.128) {$0.01$};
\draw[gp path] (1.688,2.472)--(1.778,2.472);
\draw[gp path] (1.688,2.927)--(1.778,2.927);
\draw[gp path] (1.688,3.161)--(1.778,3.161);
\draw[gp path] (1.688,3.271)--(1.868,3.271);
\node[gp node right] at (1.504,3.271) {$0.1$};
\draw[gp path] (1.688,3.616)--(1.778,3.616);
\draw[gp path] (1.688,4.070)--(1.778,4.070);
\draw[gp path] (1.688,4.304)--(1.778,4.304);
\draw[gp path] (1.688,4.415)--(1.868,4.415);
\node[gp node right] at (1.504,4.415) {$1$};
\draw[gp path] (1.688,4.759)--(1.778,4.759);
\draw[gp path] (1.688,5.214)--(1.778,5.214);
\draw[gp path] (1.688,5.447)--(1.778,5.447);
\draw[gp path] (1.688,5.558)--(1.868,5.558);
\node[gp node right] at (1.504,5.558) {$10$};
\draw[gp path] (1.688,5.902)--(1.778,5.902);
\draw[gp path] (1.688,6.357)--(1.778,6.357);
\draw[gp path] (1.688,6.590)--(1.778,6.590);
\draw[gp path] (1.688,6.701)--(1.868,6.701);
\node[gp node right] at (1.504,6.701) {$100$};
\draw[gp path] (1.688,0.985)--(1.688,1.165);
\draw[gp path] (1.688,6.701)--(1.688,6.521);
\node[gp node center] at (1.688,0.677) {$3$};
\draw[gp path] (3.200,0.985)--(3.200,1.165);
\draw[gp path] (3.200,6.701)--(3.200,6.521);
\node[gp node center] at (3.200,0.677) {$5$};
\draw[gp path] (4.712,0.985)--(4.712,1.165);
\draw[gp path] (4.712,6.701)--(4.712,6.521);
\node[gp node center] at (4.712,0.677) {$7$};
\draw[gp path] (6.223,0.985)--(6.223,1.165);
\draw[gp path] (6.223,6.701)--(6.223,6.521);
\node[gp node center] at (6.223,0.677) {$9$};
\draw[gp path] (7.735,0.985)--(7.735,1.165);
\draw[gp path] (7.735,6.701)--(7.735,6.521);
\node[gp node center] at (7.735,0.677) {$11$};
\draw[gp path] (7.735,0.985)--(7.555,0.985);
\node[gp node left] at (7.919,0.985) {$10^{0}$};
\draw[gp path] (7.735,1.557)--(7.555,1.557);
\node[gp node left] at (7.919,1.557) {$10^{1}$};
\draw[gp path] (7.735,2.128)--(7.555,2.128);
\node[gp node left] at (7.919,2.128) {$10^{2}$};
\draw[gp path] (7.735,2.700)--(7.555,2.700);
\node[gp node left] at (7.919,2.700) {$10^{3}$};
\draw[gp path] (7.735,3.271)--(7.555,3.271);
\node[gp node left] at (7.919,3.271) {$10^{4}$};
\draw[gp path] (7.735,3.843)--(7.555,3.843);
\node[gp node left] at (7.919,3.843) {$10^{5}$};
\draw[gp path] (7.735,4.415)--(7.555,4.415);
\node[gp node left] at (7.919,4.415) {$10^{6}$};
\draw[gp path] (7.735,4.986)--(7.555,4.986);
\node[gp node left] at (7.919,4.986) {$10^{7}$};
\draw[gp path] (7.735,5.558)--(7.555,5.558);
\node[gp node left] at (7.919,5.558) {$10^{8}$};
\draw[gp path] (7.735,6.129)--(7.555,6.129);
\node[gp node left] at (7.919,6.129) {$10^{9}$};
\draw[gp path] (7.735,6.701)--(7.555,6.701);
\node[gp node left] at (7.919,6.701) {$10^{10}$};
\draw[gp path] (1.688,6.701)--(1.688,0.985)--(7.735,0.985)--(7.735,6.701)--cycle;
\node[gp node center,rotate=-270] at (0.246,3.843) {$\hat{R}_{r\infty}$ and $N^*$};
\node[gp node center,rotate=-270] at (8.992,3.843) {$\hat{R}_{r0}$};
\node[gp node center] at (4.711,0.215) {$n$};
\node[gp node right] at (6.267,4.293) {$\hat{R}_{r\infty}$};
\draw[gp path] (6.451,4.293)--(7.367,4.293);
\draw[gp path] (1.688,4.101)--(3.200,3.561)--(4.712,2.983)--(6.223,2.383)--(7.735,1.767);
\gpsetpointsize{4.00}
\gppoint{gp mark 7}{(1.688,4.101)}
\gppoint{gp mark 7}{(3.200,3.561)}
\gppoint{gp mark 7}{(4.712,2.983)}
\gppoint{gp mark 7}{(6.223,2.383)}
\gppoint{gp mark 7}{(7.735,1.767)}
\gppoint{gp mark 7}{(6.909,4.293)}
\node[gp node right] at (6.267,3.843) {$N^*$};
\gpsetdashtype{gp dt 2}
\draw[gp path] (6.451,3.843)--(7.367,3.843);
\draw[gp path] (1.688,5.437)--(3.200,5.690)--(4.712,5.836)--(6.223,5.946)--(7.735,6.034);
\gppoint{gp mark 5}{(1.688,5.437)}
\gppoint{gp mark 5}{(3.200,5.690)}
\gppoint{gp mark 5}{(4.712,5.836)}
\gppoint{gp mark 5}{(6.223,5.946)}
\gppoint{gp mark 5}{(7.735,6.034)}
\gppoint{gp mark 5}{(6.909,3.843)}
\node[gp node right] at (6.267,3.393) {$\hat{R}_{r0}$};
\gpsetdashtype{gp dt 4}
\draw[gp path] (6.451,3.393)--(7.367,3.393);
\draw[gp path] (1.688,1.516)--(3.200,2.623)--(4.712,3.828)--(6.223,5.131)--(7.735,6.520);
\gppoint{gp mark 9}{(1.688,1.516)}
\gppoint{gp mark 9}{(3.200,2.623)}
\gppoint{gp mark 9}{(4.712,3.828)}
\gppoint{gp mark 9}{(6.223,5.131)}
\gppoint{gp mark 9}{(7.735,6.520)}
\gppoint{gp mark 9}{(6.909,3.393)}
\gpsetdashtype{gp dt solid}
\draw[gp path] (1.688,6.701)--(1.688,0.985)--(7.735,0.985)--(7.735,6.701)--cycle;
%% coordinates of the plot area
\gpdefrectangularnode{gp plot 1}{\pgfpoint{1.688cm}{0.985cm}}{\pgfpoint{7.735cm}{6.701cm}}
\end{tikzpicture}
%% gnuplot variables
}
    \end{minipage}
    \begin{minipage}{.235\textwidth}
      {\begin{tabular}{clr}
\toprule
$n$ &	$\hat{R}_{r\infty}$ & $N^*$	\\ \midrule
3	&   0.53  	      & 8       \\
5	&   0.18	      & 14	    \\
7	&   0.056	      & 18      \\
9	&   0.017	      & 22      \\
11	&   0.0048	      & 27      \\
\bottomrule
\end{tabular}
}
    \end{minipage}

\caption{Asymptotic and initial error ratios, and transition point for rules up to order eleven. 
 The asymptotic error ratio decreases exponentially while the initial error ratio increases also exponentially.
 The transition point increases almost linearly.}
\end{figure}

%\begin{samepage}
In \cref{fig:utschVSnewtoncotes}, we show, for rules up to eleventh order, the
asymptotic error constant ratio
\begin{equation}
    \hat{R}_{r_\infty} \triangleq \lim_{N\rightarrow\infty} 
    \left| \frac{\hat{R}_\text{new family}}{\hat{R}_\text{Newton-Cotes}} \right|
    = \left| \frac{(\hat{R}_n)_\text{new family}}{(\hat{R}_n)_\text{Newton-Cotes}} \right|
\end{equation}
as well as the initial error ratio
\begin{equation}
    \hat{R}_{r0} \triangleq
    \left|
    \frac{\hat{R}_\text{new family}}{\hat{R}_\text{Newton-Cotes}}
    \right|_{N=n}
    = \left| \frac{(\hat{R}_n)_\text{new family}}{(\hat{R}_n)_\text{Newton-Cotes}} \right|
    (n-1)^{n+1}
\end{equation}
and the transition point $N^*$, defined as 
\begin{equation}
    N^* \triangleq \min \{N \in \N | \hat{R}_\text{new family} < \hat{R}_\text{Newton-Cotes}\}
\end{equation}
i.e.\
\begin{equation}
    N^* = 
        \left\lceil
        \frac{1-(n-1)^2/(\hat{R}_{r0})^{1/(n-1)}}
             {1-(n-1)/(\hat{R}_{r0})^{1/(n-1)}}
        \right\rceil
\end{equation}
%\end{samepage}

When the rule order is increased, the asymptotic error ratio,
$\hat{R}_{r\infty}$, decreases, and it is always less than one. Using the
proposed rules is therefor advantageous when evaluating a large number of
points. In fact, the number of points at which these rules start performing
better than their Newton-Cotes pairs, $N^*$, is relatively small and
increases moderately with the rule order. In a single step, on the other
hand, the initial error ratio $C_{r0}$ is very big and shows that
Newton-Cotes rules perform much better in this scenario. For this reason, the
use of the proposed rules with less than $N^*$ points is strongly disadvised.

\subsection{Error estimation}

\begin{figure}
    \label{plot:error_estimation}
    \centering
    \input{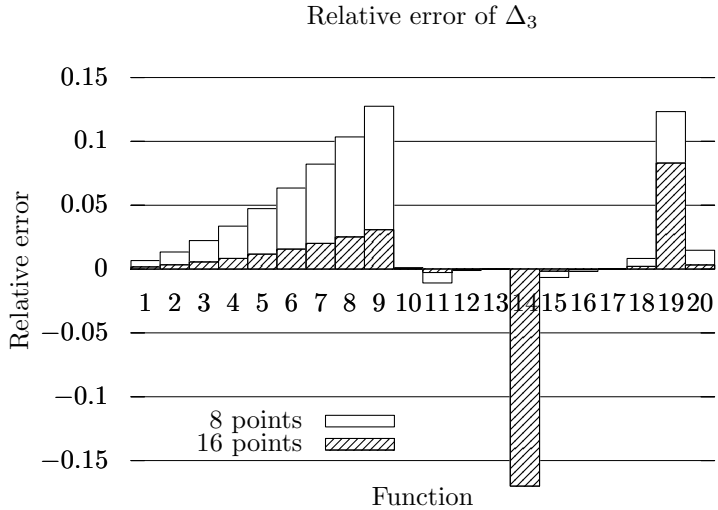}
    \caption{Performance of the correction term $\Delta_3$ as a predictor of the
    midpoint rule's error for the test functions defined and numbered in
    \cref{sec:numerical_results}. The relative error is defined as $\frac
    {\Delta_3 - R(f)}{R(f)}$, where $R(f)$ is the integration error of the
    midpoint rule.}

\end{figure}

In \cref{plot:error_estimation}, we show the performance of the correction term $\Delta_3$ as a predictor of the midpoint rule's error. It is mostly
accurate within a 15\% deviation and tends to overestimate the error. For the
\nth{3} order rule, it grossly overestimates the error. The exception is for
function 14, which is periodic and integrated in its period, meaning that
$\Delta_3$ is zero and so it loses its worthiness as an error
estimate. This is expected as low order interpolatory rules are very
accurate for periodic functions.

Functions 1 to 9 are polynomials, and the relative error rises with the degree
of the integrand, as should be expected. The relative error peaks again in
function~19, as it has a singularity on the first derivative. 
Function~20 has a singularity on the second derivative,
 but this does not have a pronounced effect on the relative error of the estimator. 
The relative errors are smaller for the other functions, which include trigonometric
and exponential functions. See \cref{sec:numerical_results} for the definition
of the test functions.

\section{Conclusion}
\label{sec:conclusion}
In this article we presented a family of uniformly spaced interpolatory
integration rules which use function evaluations in and beyond each integration
sub-interval.

These rules are based on polynomial approximations of arbitrary degree,
 and use only one function evaluation inside each integration sub-interval.
 The price paid for higher degree interpolations is one extra function evaluation outside the integration interval per degree increased.
 This allows the rules to be expressed as the midpoint rule with an additional correction term.
 This term depends on the evaluation of the integrand on a small fixed number of points close to the integration boundaries,
 regardless of the total number of function evaluations used by the associated midpoint rule.
 This correction term may also be used for quickly estimating, \emph{a priori,} the total number of function evaluations needed to comply with a global precision criterion.

This family of rules show a smaller error constant for the Peano's error formula than the corresponding constant for Newton-Cotes rules beginning from a relatively small number of function evaluations for practical applications.

As an alternative to the third order rule, we propose a modified rule 
 that uses the function values at the interval's extrema 
 instead of points outside the integration interval,
 while preserving the characteristic of ``midpoint rule with added correction term''.
Since this alternative rule only needs function evaluations inside 
 the integration interval, it rule covers more practical applications.
In the numerical results, this modified rule had smaller errors than the widely used Simpson's rule, 
but slightly higher errors than the originally proposed third order rule.
This modification can be easily extended to higher order rules.

The numerical results obtained were very satisfactory and consistent with the theoretical analysis, thus encouraging the adoption of these rules for the evaluation of integrals in practical problems.

\bibliographystyle{siamplain}
\bibliography{bib}

\pagebreak
\appendix

\section{Experimental results}
\label{sec:numerical_results}

Unless noted otherwise, the exact result of the integrals is 1.
{ \small
\[ \begin{tabular*}{\textwidth}{>{$}c<{$} @{\extracolsep{\fill}} D{.}{.}{2,8} @{\extracolsep{\fill}} D{.}{.}{2,8} @{\extracolsep{\fill}} D{.}{.}{2,8} @{\extracolsep{\fill}} D{.}{.}{2,8} @{\extracolsep{\fill}} D{.}{.}{2,8}}
\toprule 
n& \multicolumn{5}{c}{$\displaystyle \text{\tiny{(1)}} \int_{0}^{1}5x^{4}\,\dif x$} \\ 
& \multicolumn{1}{c}{\sc Midpoint} &	\multicolumn{1}{c}{\sc Simpson} &	\multicolumn{1}{c}{\sc \nth{3} Order} &	\multicolumn{1}{c}{\sc \nth{3} Order interval} &	\multicolumn{1}{c}{\sc Derivatives} \\ 
\midrule 
9 & 0.98973416 & 1.00016276 & 1.00014751 & 0.99983762 & 1.00006074\\ 
17 & 0.99711824 & 1.00001017 & 1.00000700 & 0.99999136 & 1.00000288\\ 
33 & 0.99923489 & 1.00000064 & 1.00000038 & 0.99999950 & 1.00000016\\ 
65 & 0.99980277 & 1.00000004 & 1.00000002 & 0.99999997 & 1.00000001\\ 
129 & 0.99994992 & 1.00000000 & 1.00000000 &0.10000000 & 1.00000000\\ 
\bottomrule\end{tabular*} \] 

\[ \begin{tabular*}{\textwidth}{>{$}c<{$} @{\extracolsep{\fill}} D{.}{.}{2,8} @{\extracolsep{\fill}} D{.}{.}{2,8} @{\extracolsep{\fill}} D{.}{.}{2,8} @{\extracolsep{\fill}} D{.}{.}{2,8} @{\extracolsep{\fill}} D{.}{.}{2,8}}
\toprule 
n& \multicolumn{5}{c}{$\displaystyle \text{\tiny{(2)}} \int_{0}^{1}6x^{5}\,\dif x$} \\ 
& \multicolumn{1}{c}{\sc Midpoint} &	\multicolumn{1}{c}{\sc Simpson} &	\multicolumn{1}{c}{\sc \nth{3} Order} &	\multicolumn{1}{c}{\sc \nth{3} Order interval} &	\multicolumn{1}{c}{\sc Derivatives} \\ 
\midrule 
9 & 0.98463458 & 1.00048828 & 1.00044252 & 0.99951285 & 1.00018222\\ 
17 & 0.99567998 & 1.00003052 & 1.00002099 & 0.99997407 & 1.00000864\\ 
33 & 0.99885253 & 1.00000191 & 1.00000115 & 0.99999851 & 1.00000047\\ 
65 & 0.99970417 & 1.00000012 & 1.00000007 & 0.99999991 & 1.00000003\\ 
129 & 0.99992489 & 1.00000001 & 1.00000000 & 0.99999999 & 1.00000000\\ 
\bottomrule\end{tabular*} \] 

\[ \begin{tabular*}{\textwidth}{>{$}c<{$} @{\extracolsep{\fill}} D{.}{.}{2,8} @{\extracolsep{\fill}} D{.}{.}{2,8} @{\extracolsep{\fill}} D{.}{.}{2,8} @{\extracolsep{\fill}} D{.}{.}{2,8} @{\extracolsep{\fill}} D{.}{.}{2,8}}
\toprule 
n& \multicolumn{5}{c}{$\displaystyle \text{\tiny{(3)}} \int_{0}^{1}7x^{6}\,\dif x$} \\ 
& \multicolumn{1}{c}{\sc Midpoint} &	\multicolumn{1}{c}{\sc Simpson} &	\multicolumn{1}{c}{\sc \nth{3} Order} &	\multicolumn{1}{c}{\sc \nth{3} Order interval} &	\multicolumn{1}{c}{\sc Derivatives} \\ 
\midrule 
9 & 0.97855035 & 1.00112661 & 1.00103211 & 0.99895850 & 1.00042380\\ 
17 & 0.99395685 & 1.00007101 & 1.00004897 & 0.99994201 & 1.00002015\\ 
33 & 0.99839388 & 1.00000445 & 1.00000268 & 0.99999659 & 1.00000111\\ 
65 & 0.99958586 & 1.00000028 & 1.00000016 & 0.99999979 & 1.00000006\\ 
129 & 0.99989484 & 1.00000002 & 1.00000001 & 0.99999999 & 1.00000000\\ 
\bottomrule\end{tabular*} \] 

\[ \begin{tabular*}{\textwidth}{>{$}c<{$} @{\extracolsep{\fill}} D{.}{.}{2,8} @{\extracolsep{\fill}} D{.}{.}{2,8} @{\extracolsep{\fill}} D{.}{.}{2,8} @{\extracolsep{\fill}} D{.}{.}{2,8} @{\extracolsep{\fill}} D{.}{.}{2,8}}
\toprule 
n& \multicolumn{5}{c}{$\displaystyle \text{\tiny{(4)}} \int_{0}^{1}8x^{7}\,\dif x$} \\ 
& \multicolumn{1}{c}{\sc Midpoint} &	\multicolumn{1}{c}{\sc Simpson} &	\multicolumn{1}{c}{\sc \nth{3} Order} &	\multicolumn{1}{c}{\sc \nth{3} Order interval} &	\multicolumn{1}{c}{\sc Derivatives} \\ 
\midrule 
9 & 0.97150338 & 1.00222778 & 1.00206334 & 0.99810736 & 1.00084485\\ 
17 & 0.99195060 & 1.00014162 & 1.00009792 & 0.99988903 & 1.00004027\\ 
33 & 0.99785908 & 1.00000889 & 1.00000537 & 0.99999333 & 1.00000221\\ 
65 & 0.99944785 & 1.00000056 & 1.00000031 & 0.99999959 & 1.00000013\\ 
129 & 0.99985979 & 1.00000003 & 1.00000002 & 0.99999997 & 1.00000001\\ 
\bottomrule\end{tabular*} \] 

\[ \begin{tabular*}{\textwidth}{>{$}c<{$} @{\extracolsep{\fill}} D{.}{.}{2,8} @{\extracolsep{\fill}} D{.}{.}{2,8} @{\extracolsep{\fill}} D{.}{.}{2,8} @{\extracolsep{\fill}} D{.}{.}{2,8} @{\extracolsep{\fill}} D{.}{.}{2,8}}
\toprule 
n& \multicolumn{5}{c}{$\displaystyle \text{\tiny{(5)}} \int_{0}^{1}9x^{8}\,\dif x$} \\ 
& \multicolumn{1}{c}{\sc Midpoint} &	\multicolumn{1}{c}{\sc Simpson} &	\multicolumn{1}{c}{\sc \nth{3} Order} &	\multicolumn{1}{c}{\sc \nth{3} Order interval} &	\multicolumn{1}{c}{\sc Derivatives} \\ 
\midrule 
9 & 0.96351945 & 1.00395048 & 1.00371195 & 0.99690713 & 1.00151420\\ 
17 & 0.98966330 & 1.00025397 & 1.00017624 & 0.99980889 & 1.00007242\\ 
33 & 0.99724828 & 1.00001598 & 1.00000966 & 0.99998824 & 1.00000398\\ 
65 & 0.99929015 & 1.00000100 & 1.00000057 & 0.99999927 & 1.00000023\\ 
129 & 0.99981974 & 1.00000006 & 1.00000003 & 0.99999995 & 1.00000001\\ 
\bottomrule\end{tabular*} \] 

\[ \begin{tabular*}{\textwidth}{>{$}c<{$} @{\extracolsep{\fill}} D{.}{.}{2,8} @{\extracolsep{\fill}} D{.}{.}{2,8} @{\extracolsep{\fill}} D{.}{.}{2,8} @{\extracolsep{\fill}} D{.}{.}{2,8} @{\extracolsep{\fill}} D{.}{.}{2,8}}
\toprule 
n& \multicolumn{5}{c}{$\displaystyle \text{\tiny{(6)}} \int_{0}^{1}10x^{9}\,\dif x$} \\ 
& \multicolumn{1}{c}{\sc Midpoint} &	\multicolumn{1}{c}{\sc Simpson} &	\multicolumn{1}{c}{\sc \nth{3} Order} &	\multicolumn{1}{c}{\sc \nth{3} Order interval} &	\multicolumn{1}{c}{\sc Derivatives} \\ 
\midrule 
9 & 0.95462817 & 1.00646198 & 1.00618235 & 0.99532040 & 1.00251011\\ 
17 & 0.98709736 & 1.00042131 & 1.00029369 & 0.99969525 & 1.00012056\\ 
33 & 0.99656163 & 1.00002661 & 1.00001611 & 0.99998081 & 1.00000663\\ 
65 & 0.99911277 & 1.00000167 & 1.00000094 & 0.99999880 & 1.00000039\\ 
129 & 0.99977468 & 1.00000010 & 1.00000006 & 0.99999992 & 1.00000002\\ 
\bottomrule\end{tabular*} \] 

\[ \begin{tabular*}{\textwidth}{>{$}c<{$} @{\extracolsep{\fill}} D{.}{.}{2,8} @{\extracolsep{\fill}} D{.}{.}{2,8} @{\extracolsep{\fill}} D{.}{.}{2,8} @{\extracolsep{\fill}} D{.}{.}{2,8} @{\extracolsep{\fill}} D{.}{.}{2,8}}
\toprule 
n& \multicolumn{5}{c}{$\displaystyle \text{\tiny{(7)}} \int_{0}^{1}11x^{10}\,\dif x$} \\ 
& \multicolumn{1}{c}{\sc Midpoint} &	\multicolumn{1}{c}{\sc Simpson} &	\multicolumn{1}{c}{\sc \nth{3} Order} &	\multicolumn{1}{c}{\sc \nth{3} Order interval} &	\multicolumn{1}{c}{\sc Derivatives} \\ 
\midrule 
9 & 0.94486271 & 1.00993023 & 1.00970740 & 0.99332429 & 1.00391911\\ 
17 & 0.98425551 & 1.00065838 & 1.00046143 & 0.99954178 & 1.00018919\\ 
33 & 0.99579935 & 1.00004176 & 1.00002531 & 0.99997048 & 1.00001041\\ 
65 & 0.99891573 & 1.00000262 & 1.00000148 & 0.99999813 & 1.00000061\\ 
129 & 0.99972461 & 1.00000016 & 1.00000009 & 0.99999988 & 1.00000004\\ 
\bottomrule\end{tabular*} \] 

\[ \begin{tabular*}{\textwidth}{>{$}c<{$} @{\extracolsep{\fill}} D{.}{.}{2,8} @{\extracolsep{\fill}} D{.}{.}{2,8} @{\extracolsep{\fill}} D{.}{.}{2,8} @{\extracolsep{\fill}} D{.}{.}{2,8} @{\extracolsep{\fill}} D{.}{.}{2,8}}
\toprule 
n& \multicolumn{5}{c}{$\displaystyle \text{\tiny{(8)}} \int_{0}^{1}12x^{11}\,\dif x$} \\ 
& \multicolumn{1}{c}{\sc Midpoint} &	\multicolumn{1}{c}{\sc Simpson} &	\multicolumn{1}{c}{\sc \nth{3} Order} &	\multicolumn{1}{c}{\sc \nth{3} Order interval} &	\multicolumn{1}{c}{\sc Derivatives} \\ 
\midrule 
9 & 0.93425954 & 1.01451584 & 1.01454812 & 0.99091010 & 1.00583479\\ 
17 & 0.98114084 & 1.00098118 & 1.00069200 & 0.99934231 & 1.00028332\\ 
33 & 0.99496165 & 1.00006253 & 1.00003796 & 0.99995665 & 1.00001561\\ 
65 & 0.99869903 & 1.00000393 & 1.00000223 & 0.99999723 & 1.00000092\\ 
129 & 0.99966954 & 1.00000025 & 1.00000013 & 0.99999983 & 1.00000006\\ 
\bottomrule\end{tabular*} \] 

\[ \begin{tabular*}{\textwidth}{>{$}c<{$} @{\extracolsep{\fill}} D{.}{.}{2,8} @{\extracolsep{\fill}} D{.}{.}{2,8} @{\extracolsep{\fill}} D{.}{.}{2,8} @{\extracolsep{\fill}} D{.}{.}{2,8} @{\extracolsep{\fill}} D{.}{.}{2,8}}
\toprule 
n& \multicolumn{5}{c}{$\displaystyle \text{\tiny{(9)}} \int_{0}^{1}13x^{12}\,\dif x$} \\ 
& \multicolumn{1}{c}{\sc Midpoint} &	\multicolumn{1}{c}{\sc Simpson} &	\multicolumn{1}{c}{\sc \nth{3} Order} &	\multicolumn{1}{c}{\sc \nth{3} Order interval} &	\multicolumn{1}{c}{\sc Derivatives} \\ 
\midrule 
9 & 0.92285810 & 1.02036558 & 1.02099359 & 0.98808254 & 1.00835650\\ 
17 & 0.97775670 & 1.00140680 & 1.00099931 & 0.99909096 & 1.00040848\\ 
33 & 0.99404878 & 1.00009016 & 1.00005483 & 0.99993870 & 1.00002254\\ 
65 & 0.99846271 & 1.00000567 & 1.00000321 & 0.99999604 & 1.00000132\\ 
129 & 0.99960947 & 1.00000035 & 1.00000019 & 0.99999975 & 1.00000008\\ 
\bottomrule\end{tabular*} \] 

\[ \begin{tabular*}{\textwidth}{>{$}c<{$} @{\extracolsep{\fill}} D{.}{.}{2,8} @{\extracolsep{\fill}} D{.}{.}{2,8} @{\extracolsep{\fill}} D{.}{.}{2,8} @{\extracolsep{\fill}} D{.}{.}{2,8} @{\extracolsep{\fill}} D{.}{.}{2,8}}
\toprule 
n& \multicolumn{5}{c}{$\displaystyle \text{\tiny{(10)}} \int_{0}^{1}e^x\,\dif x$} \\ 
& \multicolumn{1}{c}{\sc Midpoint} &	\multicolumn{1}{c}{\sc Simpson} &	\multicolumn{1}{c}{\sc \nth{3} Order} &	\multicolumn{1}{c}{\sc \nth{3} Order interval} &	\multicolumn{1}{c}{\sc Derivatives} \\ 
\midrule 
9 & 1.71739826 & 1.71828415 & 1.71828394 & 1.71827954 & 1.71828270\\ 
17 & 1.71803412 & 1.71828197 & 1.71828193 & 1.71828171 & 1.71828187\\ 
33 & 1.71821609 & 1.71828184 & 1.71828183 & 1.71828182 & 1.71828183\\ 
65 & 1.71826488 & 1.71828183 & 1.71828183 & 1.71828183 & 1.71828183\\ 
\midrule 
\text{Exact value} & 1.71828183 & 1.71828183 & 1.71828183 & 1.71828183 & 1.71828183\\ 
\bottomrule\end{tabular*} \] 

\[ \begin{tabular*}{\textwidth}{>{$}c<{$} @{\extracolsep{\fill}} D{.}{.}{2,8} @{\extracolsep{\fill}} D{.}{.}{2,8} @{\extracolsep{\fill}} D{.}{.}{2,8} @{\extracolsep{\fill}} D{.}{.}{2,8} @{\extracolsep{\fill}} D{.}{.}{2,8}}
\toprule 
n& \multicolumn{5}{c}{$\displaystyle \text{\tiny{(11)}} \int_{0}^{1}\sin(\pi x)\,\dif x$} \\ 
& \multicolumn{1}{c}{\sc Midpoint} &	\multicolumn{1}{c}{\sc Simpson} &	\multicolumn{1}{c}{\sc \nth{3} Order} &	\multicolumn{1}{c}{\sc \nth{3} Order interval} &	\multicolumn{1}{c}{\sc Derivatives} \\ 
\midrule 
9 & 0.63986339 & 0.63670545 & 0.63669606 & 0.63652116 & 0.63665133\\ 
17 & 0.63752656 & 0.63662505 & 0.63662339 & 0.63661493 & 0.63662126\\ 
33 & 0.63686024 & 0.63662010 & 0.63661997 & 0.63661950 & 0.63661985\\ 
65 & 0.63668174 & 0.63661979 & 0.63661978 & 0.63661976 & 0.63661978\\ 
\midrule 
\text{Exact value} & 0.63661977 & 0.63661977 & 0.63661977 & 0.63661977 & 0.63661977\\ 
\bottomrule\end{tabular*} \] 

\[ \begin{tabular*}{\textwidth}{>{$}c<{$} @{\extracolsep{\fill}} D{.}{.}{2,8} @{\extracolsep{\fill}} D{.}{.}{2,8} @{\extracolsep{\fill}} D{.}{.}{2,8} @{\extracolsep{\fill}} D{.}{.}{2,8} @{\extracolsep{\fill}} D{.}{.}{2,8}}
\toprule 
n& \multicolumn{5}{c}{$\displaystyle \text{\tiny{(12)}} \int_{0}^{1}\cos(x)\,\dif x$} \\ 
& \multicolumn{1}{c}{\sc Midpoint} &	\multicolumn{1}{c}{\sc Simpson} &	\multicolumn{1}{c}{\sc \nth{3} Order} &	\multicolumn{1}{c}{\sc \nth{3} Order interval} &	\multicolumn{1}{c}{\sc Derivatives} \\ 
\midrule 
9 & 0.84190400 & 0.84147213 & 0.84147202 & 0.84146983 & 0.84147141\\ 
17 & 0.84159232 & 0.84147106 & 0.84147103 & 0.84147092 & 0.84147101\\ 
33 & 0.84150318 & 0.84147099 & 0.84147099 & 0.84147098 & 0.84147099\\ 
65 & 0.84147928 & 0.84147099 & 0.84147098 & 0.84147098 & 0.84147098\\ 
\midrule 
\text{Exact value} & 0.84147098 & 0.84147098 & 0.84147098 & 0.84147098 & 0.84147098\\ 
\bottomrule\end{tabular*} \] 

\[ \begin{tabular*}{\textwidth}{>{$}c<{$} @{\extracolsep{\fill}} D{.}{.}{2,8} @{\extracolsep{\fill}} D{.}{.}{2,8} @{\extracolsep{\fill}} D{.}{.}{2,8} @{\extracolsep{\fill}} D{.}{.}{2,8} @{\extracolsep{\fill}} D{.}{.}{2,8}}
\toprule 
n& \multicolumn{5}{c}{$\displaystyle \text{\tiny{(13)}} \int_{0}^{1}\frac{1}{1+x^2}\,\dif x$} \\ 
& \multicolumn{1}{c}{\sc Midpoint} &	\multicolumn{1}{c}{\sc Simpson} &	\multicolumn{1}{c}{\sc \nth{3} Order} &	\multicolumn{1}{c}{\sc \nth{3} Order interval} &	\multicolumn{1}{c}{\sc Derivatives} \\ 
\midrule 
9 & 0.78565536 & 0.78539813 & 0.78539816 & 0.78540111 & 0.78539816\\ 
17 & 0.78547025 & 0.78539816 & 0.78539816 & 0.78539823 & 0.78539816\\ 
33 & 0.78541729 & 0.78539816 & 0.78539816 & 0.78539817 & 0.78539816\\ 
65 & 0.78540309 & 0.78539816 & 0.78539816 & 0.78539816 & 0.78539816\\ 
\midrule 
\text{Exact value} & 0.78539816 & 0.78539816 & 0.78539816 & 0.78539816 & 0.78539816\\ 
\bottomrule\end{tabular*} \] 

\[ \begin{tabular*}{\textwidth}{>{$}c<{$} @{\extracolsep{\fill}} D{.}{.}{2,8} @{\extracolsep{\fill}} D{.}{.}{2,8} @{\extracolsep{\fill}} D{.}{.}{2,8} @{\extracolsep{\fill}} D{.}{.}{2,8} @{\extracolsep{\fill}} D{.}{.}{2,8}}
\toprule 
n& \multicolumn{5}{c}{$\displaystyle \text{\tiny{(14)}} \int_{0}^{1}\frac{2}{2+\sin(10 \pi x)}\,\dif x$} \\ 
& \multicolumn{1}{c}{\sc Midpoint} &	\multicolumn{1}{c}{\sc Simpson} &	\multicolumn{1}{c}{\sc \nth{3} Order} &	\multicolumn{1}{c}{\sc \nth{3} Order interval} &	\multicolumn{1}{c}{\sc Derivatives} \\ 
\midrule 
9 & 1.15470054 & 1.15079365 & 1.15470052 & 1.14845436 & 1.15470052\\ 
17 & 1.15470054 & 1.15468009 & 1.15384615 & 1.15000000 & 1.15384615\\ 
33 & 1.15470054 & 1.15470054 & 1.15470054 & 1.15449366 & 1.15470054\\ 
65 & 1.15470054 & 1.15470054 & 1.15470054 & 1.15469686 & 1.15470054\\ 
\midrule 
\text{Exact value} & 1.15470054 & 1.15470054 & 1.15470054 & 1.15470054 & 1.15470054\\ 
\bottomrule\end{tabular*} \] 

\[ \begin{tabular*}{\textwidth}{>{$}c<{$} @{\extracolsep{\fill}} D{.}{.}{2,8} @{\extracolsep{\fill}} D{.}{.}{2,8} @{\extracolsep{\fill}} D{.}{.}{2,8} @{\extracolsep{\fill}} D{.}{.}{2,8} @{\extracolsep{\fill}} D{.}{.}{2,8}}
\toprule 
n& \multicolumn{5}{c}{$\displaystyle \text{\tiny{(15)}} \int_{0}^{1}\frac{1}{1+x^4}\,\dif x$} \\ 
& \multicolumn{1}{c}{\sc Midpoint} &	\multicolumn{1}{c}{\sc Simpson} &	\multicolumn{1}{c}{\sc \nth{3} Order} &	\multicolumn{1}{c}{\sc \nth{3} Order interval} &	\multicolumn{1}{c}{\sc Derivatives} \\ 
\midrule 
9 & 0.86748850 & 0.86698105 & 0.86698036 & 0.86695345 & 0.86697602\\ 
17 & 0.86711725 & 0.86697350 & 0.86697334 & 0.86697229 & 0.86697313\\ 
33 & 0.86701125 & 0.86697302 & 0.86697301 & 0.86697296 & 0.86697300\\ 
65 & 0.86698285 & 0.86697299 & 0.86697299 & 0.86697299 & 0.86697299\\ 
\midrule 
\text{Exact value} & 0.86697299 & 0.86697299 & 0.86697299 & 0.86697299 & 0.86697299\\ 
\bottomrule\end{tabular*} \] 

\[ \begin{tabular*}{\textwidth}{>{$}c<{$} @{\extracolsep{\fill}} D{.}{.}{2,8} @{\extracolsep{\fill}} D{.}{.}{2,8} @{\extracolsep{\fill}} D{.}{.}{2,8} @{\extracolsep{\fill}} D{.}{.}{2,8} @{\extracolsep{\fill}} D{.}{.}{2,8}}
\toprule 
n& \multicolumn{5}{c}{$\displaystyle \text{\tiny{(16)}} \int_{0}^{1}\frac{1}{1+e^x}\,\dif x$} \\ 
& \multicolumn{1}{c}{\sc Midpoint} &	\multicolumn{1}{c}{\sc Simpson} &	\multicolumn{1}{c}{\sc \nth{3} Order} &	\multicolumn{1}{c}{\sc \nth{3} Order interval} &	\multicolumn{1}{c}{\sc Derivatives} \\ 
\midrule 
9 & 0.37985801 & 0.37988537 & 0.37988538 & 0.37988562 & 0.37988545\\ 
17 & 0.37987779 & 0.37988549 & 0.37988549 & 0.37988550 & 0.37988549\\ 
33 & 0.37988345 & 0.37988549 & 0.37988549 & 0.37988549 & 0.37988549\\ 
65 & 0.37988497 & 0.37988549 & 0.37988549 & 0.37988549 & 0.37988549\\ 
\midrule 
\text{Exact value} & 0.37988549 & 0.37988549 & 0.37988549 & 0.37988549 & 0.37988549\\ 
\bottomrule\end{tabular*} \] 

\[ \begin{tabular*}{\textwidth}{>{$}c<{$} @{\extracolsep{\fill}} D{.}{.}{2,8} @{\extracolsep{\fill}} D{.}{.}{2,8} @{\extracolsep{\fill}} D{.}{.}{2,8} @{\extracolsep{\fill}} D{.}{.}{2,8} @{\extracolsep{\fill}} D{.}{.}{2,8}}
\toprule 
n& \multicolumn{5}{c}{$\displaystyle \text{\tiny{(17)}} \int_{0}^{1}\frac{23}{25} \cosh(x) - \cos(x)\,\dif x$} \\ 
& \multicolumn{1}{c}{\sc Midpoint} &	\multicolumn{1}{c}{\sc Simpson} &	\multicolumn{1}{c}{\sc \nth{3} Order} &	\multicolumn{1}{c}{\sc \nth{3} Order interval} &	\multicolumn{1}{c}{\sc Derivatives} \\ 
\midrule 
9 & 0.23872514 & 0.23971443 & 0.23971441 & 0.23971383 & 0.23971423\\ 
17 & 0.23943692 & 0.23971413 & 0.23971413 & 0.23971410 & 0.23971412\\ 
33 & 0.23964055 & 0.23971411 & 0.23971411 & 0.23971411 & 0.23971411\\ 
65 & 0.23969515 & 0.23971411 & 0.23971411 & 0.23971411 & 0.23971411\\ 
\midrule 
\text{Exact value} & 0.23971411 & 0.23971411 & 0.23971411 & 0.23971411 & 0.23971411\\ 
\bottomrule\end{tabular*} \] 

\[ \begin{tabular*}{\textwidth}{>{$}c<{$} @{\extracolsep{\fill}} D{.}{.}{2,8} @{\extracolsep{\fill}} D{.}{.}{2,8} @{\extracolsep{\fill}} D{.}{.}{2,8} @{\extracolsep{\fill}} D{.}{.}{2,8} @{\extracolsep{\fill}} D{.}{.}{2,8}}
\toprule 
n& \multicolumn{5}{c}{$\displaystyle \text{\tiny{(18)}} \int_{0}^{1}\frac{1}{1+x}\,\dif x$} \\ 
& \multicolumn{1}{c}{\sc Midpoint} &	\multicolumn{1}{c}{\sc Simpson} &	\multicolumn{1}{c}{\sc \nth{3} Order} &	\multicolumn{1}{c}{\sc \nth{3} Order interval} &	\multicolumn{1}{c}{\sc Derivatives} \\ 
\midrule 
9 & 0.69276241 & 0.69315453 & 0.69315409 & 0.69314094 & 0.69315000\\ 
17 & 0.69303913 & 0.69314765 & 0.69314751 & 0.69314681 & 0.69314732\\ 
33 & 0.69311849 & 0.69314721 & 0.69314720 & 0.69314716 & 0.69314719\\ 
65 & 0.69313978 & 0.69314718 & 0.69314718 & 0.69314718 & 0.69314718\\ 
\midrule 
\text{Exact value} & 0.69314718 & 0.69314718 & 0.69314718 & 0.69314718 & 0.69314718\\ 
\bottomrule\end{tabular*} \] 

\[ \begin{tabular*}{\textwidth}{>{$}c<{$} @{\extracolsep{\fill}} D{.}{.}{2,8} @{\extracolsep{\fill}} D{.}{.}{2,8} @{\extracolsep{\fill}} D{.}{.}{2,8} @{\extracolsep{\fill}} D{.}{.}{2,8} @{\extracolsep{\fill}} D{.}{.}{2,8}}
\toprule 
n& \multicolumn{5}{c}{$\displaystyle \text{\tiny{(19)}} \int_{0}^{1}\sqrt{|x^2-0.25|^3}\,\dif x$} \\ 
& \multicolumn{1}{c}{\sc Midpoint} &	\multicolumn{1}{c}{\sc Simpson} &	\multicolumn{1}{c}{\sc \nth{3} Order} &	\multicolumn{1}{c}{\sc \nth{3} Order interval} &	\multicolumn{1}{c}{\sc Derivatives} \\ 
\midrule 
9 & 0.14732665 & 0.14903320 & 0.14848608 & 0.14847253 & 0.14848191\\ 
17 & 0.14845436 & 0.14889952 & 0.14881348 & 0.14881276 & 0.14881329\\ 
33 & 0.14876407 & 0.14887651 & 0.14886211 & 0.14886207 & 0.14886210\\ 
65 & 0.14884450 & 0.14887248 & 0.14887000 & 0.14887000 & 0.14887000\\ 
\midrule 
\text{Exact value} & 0.14887162 & 0.14887162 & 0.14887162 & 0.14887162 & 0.14887162\\ 
\bottomrule\end{tabular*} \] 

\[ \begin{tabular*}{\textwidth}{>{$}c<{$} @{\extracolsep{\fill}} D{.}{.}{2,8} @{\extracolsep{\fill}} D{.}{.}{2,8} @{\extracolsep{\fill}} D{.}{.}{2,8} @{\extracolsep{\fill}} D{.}{.}{2,8} @{\extracolsep{\fill}} D{.}{.}{2,8}}
\toprule 
n& \multicolumn{5}{c}{$\displaystyle \text{\tiny{(20)}} \int_{0}^{1}\sqrt{|x^2-0.25|^5}\,\dif x$} \\ 
& \multicolumn{1}{c}{\sc Midpoint} &	\multicolumn{1}{c}{\sc Simpson} &	\multicolumn{1}{c}{\sc \nth{3} Order} &	\multicolumn{1}{c}{\sc \nth{3} Order interval} &	\multicolumn{1}{c}{\sc Derivatives} \\ 
\midrule 
9 & 0.06386234 & 0.06556386 & 0.06560246 & 0.06546068 & 0.06556174\\ 
17 & 0.06504820 & 0.06551704 & 0.06551935 & 0.06551211 & 0.06551742\\ 
33 & 0.06539065 & 0.06551484 & 0.06551505 & 0.06551464 & 0.06551495\\ 
65 & 0.06548275 & 0.06551477 & 0.06551479 & 0.06551476 & 0.06551478\\ 
\midrule 
\text{Exact value} & 0.06551477 & 0.06551477 & 0.06551477 & 0.06551477 & 0.06551477\\ 
\bottomrule\end{tabular*} \] 

}

\end{document}